\documentclass[a4paper]{amsart}

\usepackage{ifthen}
\usepackage{soul}

\newcommand{\correct}[2]{{\ifthenelse{\boolean{printNew}}{{\ifthenelse{\boolean{printOld}}{\st{#1}\color{red} #2}{#2}}}{\ifthenelse{\boolean{printOld}}{#1}{#2}}}}

\provideboolean{printOld}
\provideboolean{printNew}
\setboolean{printNew}{true}
 \setboolean{printOld}{false}

\usepackage{graphics}
\usepackage{color}
\usepackage{amssymb}
\usepackage[all]{xy}
\usepackage{amsmath}
\usepackage{stmaryrd}
\definecolor{nenhanced}{rgb}{0.2,0.30,0.13}
\newcommand{\enh}[1]{{}}
\newcommand{\nenh}[1]{{#1}}

\newcommand{\kbb}{{\kss \bullet}}
\newcommand{\kst}{\,|\;}

\newcommand{\kss}{\scriptscriptstyle}

\def\interior{\operatorname{int}}
\newcommand{\toric}[1]{\mathbb T \mathbb V({#1})} % Funktor "torische Varietaet"

\newcommand{\ku}{w}  % Lokalisierungs-Grad

\CompileMatrices

\newtheorem{theorem}{Theorem}[section]
\newtheorem{lemma}[theorem]{Lemma}
\newtheorem{proposition}[theorem]{Proposition}
\newtheorem{corollary}[theorem]{Corollary}
\theoremstyle{definition}
\newtheorem{definition}[theorem]{Definition}
\newtheorem{example}[theorem]{Example}

\newtheorem{remark}[theorem]{Remark}
\theoremstyle{remark}

\def\div{{\rm div}}

\def\mal{\! \cdot \!}

\def\t#1{\widetilde{#1}}
\def\b#1{\overline{#1}}
\def\bangle#1{\langle #1 \rangle}

\def\AA{{\mathbb A}}
\def\CC{{\mathbb C}}
\def\KK{{\mathbb K}}

\def\ZZ{{\mathbb Z}}
\def\RR{{\mathbb R}}
\def\NN{{\mathbb N}}
\def\QQ{{\mathbb Q}}
\def\PP{{\mathbb P}}

\def\Of{{\mathcal{O}}}

\def\CDiv{\operatorname{CDiv}}
\def\WDiv{\operatorname{WDiv}}

\def\id{{\rm id}}

\def\Hom{{\rm Hom}}

\def\Supp{{\rm Supp}}
\def\Spec{{\rm Spec}}

\def\PPDiv{{\rm PPDiv}}
\def\EnPPDiv{{\rm EnPPDiv}}
\def\Pol{{\rm Pol}}
\def\EnPol{{\rm EnPol}}
\def\Irr{{\rm Irr}}
\def\Loc{{\rm Loc}}

\def\conv{{\rm conv}}

\def\eval{{\rm eval}}

\def\tail{{\rm tail}}

\def\face{{\rm face}}

\DeclareMathOperator{\supp}{supp}
\newcommand{\CS}{\mathcal{S}}
\newcommand{\CE}{\mathcal{E}}
\newcommand{\D}{\mathcal{D}}
\newcommand{\sigmav}{\sigma^\vee}
\newcommand{\mul}[2]{\langle #1,#2 \rangle}

\newcommand{\ppv}[1]{{X}({#1})}  % Funktor: [pp-fan] -> Varietaet

\newcounter{itemnumber}

\begin{document}

\title[Algebraic torus actions]%
{Gluing affine torus actions via divisorial fans}
\author[K.~Altmann]{Klaus Altmann}
\address{Fachbereich Mathematik und Informatik, 
Freie Universit\"at Berlin,
Arnimalle 3, 
14195 Berlin, 
Germany}
\email{altmann@math.fu-berlin.de}
\author[J.~Hausen]{J\"urgen Hausen} 
\address{Mathematisches Institut, Universit\"at T\"ubingen,
Auf der Morgenstelle 10, 72076 T\"ubingen, Germany}
\email{hausen@mail.mathematik.uni-tuebingen.de}
\author[H.~S\"uss]{Hendrik S\"uss }
\address{Institut f\"ur Mathematik,
        LS Algebra und Geometrie,
        Brandenburgische Technische Universit\"at Cottbus,
        PF 10 13 44, 
        03013 Cottbus, Germany}
\email{suess@math.tu-cottbus.de}

\begin{abstract}
Generalizing the passage from a 
fan to a toric variety, 
we provide a combinatorial approach to
construct arbitrary effective torus 
actions on normal, algebraic varieties.
Based on the notion of a 
``proper polyhedral divisor''
introduced in earlier work, we 
develop the concept 
of a ``divisorial fan'' 
and show that these objects 
encode the equivariant gluing of affine 
varieties with torus action.
We characterize separateness and 
completeness of the resulting 
varieties in terms of divisorial fans,
and we study examples like 
$\CC^*$-surfaces and projectivizations
of (non-split) vector bundles over toric     
varieties.
\end{abstract}

\maketitle

\section{Introduction}

This paper continues work of the first two 
authors~\cite{AlHa}, where the concept of 
``proper polyhedral divisors (pp-divisors)'' 
was introduced in order to provide a complete 
description of normal affine varieties $X$ 
that come with an effective action of an 
algebraic torus $T$.
Recall that such a pp-divisor lives on 
a normal semiprojective 
(e.g., affine or projective)
variety $Y$, and, at first glance,  
is just a finite linear combination
\begin{eqnarray*}
\mathcal{D}
& = & 
\sum_D \Delta_D \otimes D,
\end{eqnarray*}
where $D$ runs over the prime divisors 
of $Y$ and the coefficients $\Delta_D$ 
are convex polyhedra, all living in a 
common rational vector space $N_\QQ$ 
and all having the same pointed cone
$\sigma \subseteq N_\QQ$ as their
tail.
% Then come a couple of 
% further defining conditions on $\mathcal{D}$;
% which we will skip for the moment.
To see an example, 
let $Y$ be the projective line,
and take the points 
$0,1$ and $\infty$ as prime divisors
on $Y$. Then one obtains 
a pp-divisor $\mathcal{D}$ on $Y$ 
by prescribing polyhedral 
coefficients as follows.

\bigskip

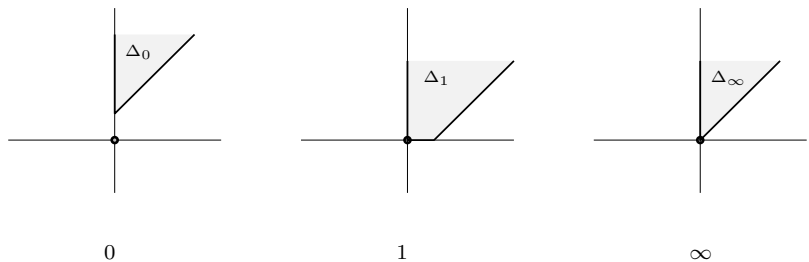
\begin{center}
\input{ppdivex.pstex_t}
\end{center}

The affine $T$-variety $X$ associated to 
$\mathcal{D}$ is the spectrum of a 
multigraded algebra~$A$ arising
from $\mathcal{D}$.
Evaluating the polyhedral coefficients
turns the pp-divisor into a piecewise 
linear map from the dual cone  
$\sigma^\vee \subseteq M_\QQ$ of 
the common tail 
to the rational Cartier divisors on $Y$: 
it sends $u \in \sigma^\vee$ to the divisor 
$\mathcal{D}(u) = \sum \alpha_D D$, 
where $\alpha_D = \min \bangle{u,\Delta_D}$.
The global sections of these evaluations 
fit together to the desired multigraded 
algebra:
\begin{eqnarray*}
A
& := & 
\bigoplus_{u \in M \cap \sigma^\vee}
\Gamma(Y,\mathcal{D}(u)).
\end{eqnarray*}

In the present paper, we pass from the affine 
case to the general one.
In the setting of toric varieties,
the general case is obtained from the affine 
one by gluing cones to a fan.
This is also our approach;
applying Sumihiro's Theorem,
we glue pp-divisors to a ``divisorial fan''.
There is an immediate naive idea of
how such a divisorial fan should look: 
all its divisors $\mathcal{D}^i$
live on the same semiprojective variety
$Y$, their polyhedral coefficients $\Delta_D^i$
live in the same vector space $N_\QQ$, and,
for every prime divisor $D$, the $\Delta_D^i$ 
should form a polyhedral subdivision. 
For example, the single pp-divisor on 
$Y = \PP^1$ discussed before could 
fit as a $\mathcal{D}^1$ 
into a divisorial fan comprising five further
pp-divisors as indicated below.

\bigskip

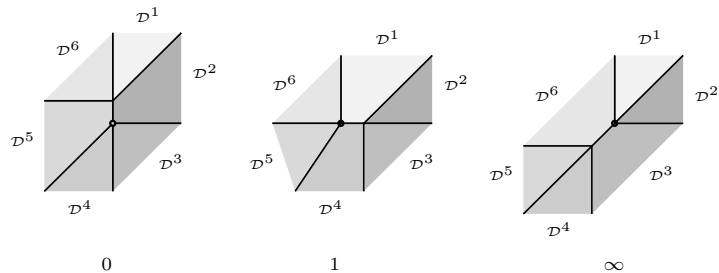
\begin{center}
\input{tangp2.pstex_t}%
\label{intro-figure}%
\end{center}

To describe the gluing of affine $T$-varieties 
amounts to understanding their open subsets in 
terms of pp-divisors; a detailed study 
is given in Sections~\ref{sec:opemb1} 
and~\ref{sec:opemb2}.
Based on this, in Section~\ref{glue} we introduce
a concept of a divisorial fan.
We show that each such divisorial fan 
canonically defines a 
normal variety with torus action (Theorem~\ref{sec:thm-divisorial-fans1}),
and it turns out that every normal variety with
effective torus action can be obtained in this way 
(Theorem~\ref{sec:thm-divisorial-fans2}).
In Section~\ref{sec:coherent}, we 
discuss ``coherent'' divisorial fans--a special
concept, 
which is much closer to the  
intuition than the general one. 
For example, the figure just drawn fits
into this framework: it describes the 
projectivization of the cotangent bundle
over the projective plane.

Following the philosophy of 
toric geometry that geometric
properties of a toric variety should 
be read off from its defining 
combinatorial data,
in Section~\ref{sec:sepcomp} we study
separateness and completeness 
and provide a complete characterization of these properties 
in terms of divisorial fans (Theorem~\ref{thm:sepcompchar}).
The last section is devoted to examples.
We give the divisorial fans of Danilov-Gizatullin 
compactifications of affine $\KK^*$-surfaces, 
recently discussed using different methods 
by Flenner, Kaliman and Zaidenberg. 
This example indicates that our constructions may be 
used for finding compactification of varieties with torus actions in a 
rather intuitive way.
In our last example we provide a translation of Klyachko's description 
of vector bundles on toric varieties 
into the picture of divisorial fans.

We expect applications of our constructions in all fields where 
toric varieties have proved their
usefulness. Following the toric program, 
next steps will be the description of divisors, bundles and equivariant maps as well as 
the understanding of intersection products and cohomology on divisorial fans.

\tableofcontents

\section{The affine case}
\label{sec:affcase}

In this section, we briefly recall the
basic concepts and results from~\cite{AlHa},
and we introduce some notions needed later.
We begin with fixing our notation 
in convex geometry.

Throughout this paper, $N$ 
denotes a lattice, i.e. a finitely generated 
free abelian group, and $M := \Hom(N,\ZZ)$ 
is the associated dual lattice.
The rational vector space associated to 
$N$ is $N_\QQ := \QQ \otimes_\ZZ N$.
Given a homomorphism $F \colon N \to N'$ 
of lattices, we write 
$F \colon  N_\QQ\to N'_\QQ$
for the corresponding linear map.
For two convex polyhedra 
$\Delta, \Delta' \subseteq N_\QQ$,
we write $\Delta \preceq \Delta'$ if 
$\Delta$ is a face of $\Delta'$.

Let $\sigma \subseteq N_\QQ$ be a 
pointed, convex, polyhedral 
cone.
A $\sigma$-polyhedron is a convex 
polyhedron
$\Delta \subseteq N_\QQ$ having
$\sigma$ as its tail cone
(also called recession cone).
With respect to Minkowski addition, 
the set  $\Pol_{\sigma}^+(N)$ of all 
$\sigma$-polyhedra is a semigroup
with cancellation law; we write 
$\Pol_{\sigma}(N)$ for the associated 
Grothendieck group. 

Now, let $Y$ be a normal, algebraic variety
defined over an algebraically closed field $\KK$. 
Except in Section \ref{sec:sepcomp}, we understand 
points always to be closed points of $Y$.
% of characteristic zero.
% (The latter seems to be not neccessary, but it was assumed in
% \cite{AlHa}).
The group of {\em polyhedral divisors\/} 
on $Y$ is defined to be
$$
\WDiv_{\QQ}(Y,\sigma)
\ := \
 \Pol_{\sigma}(N) \otimes_{\ZZ} \WDiv(Y),
$$ 
where $\WDiv_{\QQ}(Y)$ denotes the group
of rational Weil divisors on $Y$.
Any polyhedral divisor
$\mathcal{D} = \sum \Delta_D \otimes D$ 
with $\Delta_D\in\Pol^+_{\sigma}(N)$
evaluates to a piecewise linear convex 
map on the dual cone $\sigma^\vee \subseteq M_{\QQ}$
of $\sigma \subseteq N_{\QQ}$, namely
$$
\mathcal{D} \colon \sigma^\vee \to \WDiv_{\QQ}(Y),
\quad
u \mapsto \sum \eval_u(\Delta_D) D,
\text{ where }
\eval_u(\Delta_D) \; := \; \min_{v \in \Delta_D} \langle u, v\rangle.
$$
Here, convexity is understood in 
the setting of divisors, 
this means that we always
have  
$\mathcal{D}(u+u') \ge \mathcal{D}(u)+\mathcal{D}(u')$.
A {\em proper polyhedral divisor\/}
(abbreviated pp-divisor)
is a polyhedral divisor 
$\mathcal{D} \in \WDiv_{\QQ}(Y,\sigma)$
such that
\begin{enumerate}
\item 
there is a representation
$\mathcal{D} = \sum  \Delta_D\otimes D$ with
effective divisors $D \in \WDiv_\QQ(Y)$ and 
$\Delta_D \in \Pol^{+}_{\sigma}(N)$,
\item 
each evaluation $\mathcal{D}(u)$, where $u \in \sigma^\vee$,
is a semiample $\QQ$-Cartier divisor, 
i.e. has a base point free multiple,
\item 
for any $u$ in the relative interior of $\sigma^\vee$,
some multiple of $\mathcal{D}(u)$ 
is a big divisor, i.e. admits a section
with affine complement.
\end{enumerate}

>From now on, we suppose that $Y$ is, 
additionally, semiprojective,
i.e. projective over some affine variety.
Every pp-divisor 
$\mathcal{D} = \sum \Delta_D\otimes D$ 
on~$Y$ defines a sheaf
of graded $\mathcal{O}_Y$-algebras, 
and we have the 
corresponding relative spectrum: 
$$
\mathcal{A} 
\ := \
\bigoplus_{u \in M \cap \sigma^\vee}
\mathcal{O}(\mathcal{D}(u)),
\qquad
\t{X} 
\ := \
\Spec_Y(\mathcal{A}).
$$
The grading of $\mathcal{A}$ gives rise to 
an effective action of the torus $T := \Spec(\KK[M])$
on~$\t{X}$, the canonical map 
$\pi \colon \t{X} \to Y$ is a good quotient 
for this action,
and for the field
of invariant rational functions,
we have 
\begin{eqnarray*}
\KK(\t{X})^T & = & \KK(Y).
\end{eqnarray*}

By~\cite[Theorem~3.1]{AlHa},
the ring of global sections
$A := \Gamma(\t{X},\mathcal{O}) = \Gamma(Y, \mathcal{A})$
is finitely generated and normal, 
and there is a $T$-equivariant, projective, birational
morphism $r \colon \t{X} \to X$
onto the normal, affine $T$-variety
$X := X(\mathcal{D}) := \Spec(A)$.
Conversely, ~\cite[Theorem~3.4]{AlHa}
shows that every normal, affine variety
with an effective torus action arises 
in this way.

The assignment from pp-divisors to 
normal, affine varieties with torus action 
is even functorial, see~\cite[Sec.~8]{AlHa}.
Consider two pp-divisors,
$$
\mathcal{D}' 
\ = \
\sum \Delta'_{D'} \otimes D' 
\ \in \
\PPDiv_\QQ(Y',\sigma'),
\quad 
\mathcal{D} 
\ = \
\sum \Delta_D \otimes D 
\ \in \
\PPDiv_\QQ(Y,\sigma).
$$
If $\psi \colon Y' \to Y$ 
is a morphism such that
none of the supports of the 
$D$'s contains $\psi(Y)$,
and if $F \colon  N' \to N$ is 
a linear map with
$F(\sigma') \subseteq \sigma$,
then we set
$$ 
\psi^*(\mathcal{D}) 
\ := \ 
\sum_D \Delta_D\otimes \psi^\ast(D),
\qquad
F_*(\mathcal{D}') 
\ := \
\sum_{D'} \big(F(\Delta'_{D'})+\sigma\big) \otimes D'.
$$
Suppose that for some ``polyhedral principal divisor''
$\div(\mathfrak{f}) = \sum (v_i + \sigma) \otimes \div(f_i)$ 
with $v_i\in N$ and $f_i\in \KK(Y')$, 
we have inside $\WDiv_\QQ(Y',\sigma)$ the relation
$$ 
\psi^\ast(\mathcal{D})
\ \leq \
F_\ast(\mathcal{D}') 
+ 
\div(\mathfrak{f})
$$
after evaluating with arbitrary $u\in\sigma^\vee$.
Then the triple
$(\psi,F,\mathfrak{f})$
is called a map from the pp-divisor
$\mathcal{D}'$ to the pp-divisor 
$\mathcal{D}$.
It induces homomorphisms of 
$\mathcal{O}_{Y}$-modules:
\begin{equation*}
\mathcal{O}(\mathcal{D}(u))
\; \to \;
\psi_*\mathcal{O}(\mathcal{D}'(F^{*} u)),
\qquad
h 
\; \mapsto \;
\mathfrak{f}(u) \psi^*(h).
\end{equation*}
These maps fit together to a graded homomorphism
$\mathcal{A} \to \psi_*\mathcal{A}'$. 
This in turn gives rise to a commutative diagram
of equivariant morphisms,
where the rows contain the geometric 
data associated
to the pp-divisors 
$\mathcal{D}$ and $\mathcal{D}'$
respectively:
$$ 
\xymatrix{
Y 
&
{\t{X}} \ar[l]_{\pi} \ar[r]^{r}
&
X
\\
Y' \ar[u]^{\psi}
&
{\t{X}'} \ar[l]^{\pi'} \ar[r]_{r'} \ar[u]^{\t{\varphi}} 
&
X' \ar[u]_{\varphi} .
}
$$
In particular, the map $\mathcal{D}' \to \mathcal{D}$
defines an equivariant morphism $X' \to X$
with respect to $T' \to T$ defined by $F \colon N' \to N$.

Here, we will frequently consider 
a special case of the above one.
Namely, suppose that $Y'=Y$ and $\Delta'_D \subseteq \Delta_D$
holds for every prime divisor 
$D \in \WDiv(Y)$. 
Then 
$\sigma^\vee \subseteq (\sigma')^\vee$
holds for the dualized tail cones.
Moreover, for every $u \in \sigma^\vee$,
we obtain
$$ 
\mathcal{D}(u)
\ = \ 
\sum \min_{v \in \Delta_D} \langle u, v \rangle D
\ \le \
\sum \min_{v \in \Delta'_D}  \langle u, v \rangle D
\ = \ 
\mathcal{D}'(u).
$$
Consequently, we have a graded
inclusion morphism 
$\mathcal{A} \hookrightarrow \mathcal{A}'$ 
of the associated sheaves of 
$\mathcal{O}_Y$-algebras, and 
hence a monomorphism $A \hookrightarrow A'$
on the level of global sections,
which in turn determines a
$T$-equivariant 
morphism $X' \to X$.

In~\cite[Prop.~7.8 and Cor.~7.9]{AlHa},
we took a closer look at the fibers 
of the map $\pi \colon \t{X} \to Y$
arising from a pp-divisor 
$\mathcal{D} = \sum \Delta_D \otimes D$.
Suppose that all $D$'s are prime.
For a point $y \in Y$, its 
fiber polyhedron is 
the Minkowski sum
$$
\Delta_y 
\; := \;
\sum_{y \in D} \Delta_D
\; \in \; 
\Pol_{\sigma}^+(N).
$$
Let $\Lambda_y$ denote the 
normal fan of the fiber
polyhedron $\Delta_y$.
Then $\Lambda_y$ subdivides 
the cone $\sigma^\vee$, and
the faces of $\Delta_y$
are in order reversing bijection
to the cones of $\Lambda_y$ via
\begin{eqnarray*}
F 
& \mapsto &
\lambda(F)
\; := \; 
\{u \in M_{\QQ}; \; 
\langle u, v - v' \rangle \ge 0
\text{ for all } 
v \in \Delta, \, v' \in F\}.
\end{eqnarray*}
Now, for $z \in \pi^{-1}(y)$,
let $\omega(z)$ denote its orbit cone,
i.e. the convex cone generated by 
all weights $u \in M$ admitting a 
$u$-homogeneous function 
on $\pi^{-1}(y)$ with $f(u) \ne 0$.
Then there is a bijection:
$$
\{ 
T \text{-orbits in } \pi^{-1}(y)
\}
\;  \to \;
\Lambda_y
\qquad
T \mal \t{x}
\; \mapsto \;
\omega(\t{x}).
$$
This does eventually provide an order 
and dimension preserving bijection
between the $T$-orbits of $\pi^{-1}(y)$ 
and the faces of $\Delta_y$.

\enh{%
For the gluing of pp-divisors performed
later, we will broaden our notation. We
enhance the coefficient group
$\Pol_\sigma(N)$ by the empty polyhedron.

\begin{definition}
\label{sec:def-empty-set}
Let $N$ be a lattice and $\sigma \subseteq N_\QQ$ 
a pointed polyhedral cone.
Then the 
{\em enhanced group of $\sigma$-polyhedra\/}
is the set
\begin{eqnarray*}
\EnPol_\sigma(N)
& := &
\Pol_\sigma(N)
\ \cup \ 
\{ \emptyset \}
\end{eqnarray*}
together with the usual addition
and the rules
$\emptyset + \Delta := \emptyset$
and $0 \cdot \emptyset := \sigma$.
A {\em polyhedral divisor with enhanced\/}
coefficients on a normal variety $Y$ 
is a finite formal sum
$$ 
\mathcal{D} 
\ = \ 
\sum \Delta_D \otimes D,
\qquad
\text{where }
\Delta_D \in \EnPol_\sigma(N),
\quad 
D \in \WDiv_\QQ(Y) 
\text{ prime}.
$$
For a polyhedral divisor $\mathcal{D}$ 
with enhanced coefficients on $Y$,
we define its {\em irrelevant set\/}
and its {\em locus\/}
to be the subsets 
$$
\Irr(\mathcal{D}) 
\ := \ 
\bigcup_{\Delta_D = \emptyset} D
\ \subseteq \
Y,
\qquad
\Loc(\mathcal{D})
\ := \ 
Y \setminus \Irr(\mathcal{D})
\ \subseteq \
Y.
$$
We say that a polyhedral divisor $\mathcal{D}$ 
with enhanced coefficients on~$Y$ is 
a {\em pp-divisor\/} if $Y$ is semiprojective,
$\Irr(\mathcal{D})$ is the support of
an effective, semiample divisor (making
$\Loc(\mathcal{D})$ semiprojective),
and $\mathcal{D}_{\vert \Loc(\mathcal{D})}$
is a pp-divisor in the usual sense.
\end{definition}

Note that, as before, we may associate 
to any 
pp-divisor $\mathcal{D}$ with enhanced 
coefficients on a normal, semiprojective variety $Y$ 
a sheaf of graded $\mathcal{O}_Y$-algebras
\begin{eqnarray*}
\mathcal{A}
& := & 
\imath_\ast 
\left(\bigoplus_{u\in M \cap \sigma^\vee} \mathcal{O}_V
(\mathcal{D}_{|V}(u))\right),
\end{eqnarray*}
where $V = \Loc(\mathcal{D})$,
and $\imath \colon V \to Y$ is the 
inclusion.
Moreover, we have
$A := \Gamma(Y,\mathcal{A})$,
and the geometric data 
$\t{X} := \Spec_Y(\mathcal{A})$
and $X := X(\mathcal{D}) := \Spec(A)$.

Finally, we need to introduce another 
useful notion 
concerning the evaluation of the 
coefficients of a polyhedral divisor.
}%\enh

\nenh{%
Now, for the gluing of pp-divisors performed
later, it is necessary to relax our notation: We
will allow $\emptyset$ as an element of
$\Pol_\sigma(N)$.
This new element is subject to the rules $\emptyset + \Delta := \emptyset$
and $0 \cdot \emptyset := \sigma$.
Moreover,
if $\emptyset$ occurs as a coefficient of a pp-divisor
$
\,\mathcal{D} 
=
\sum \Delta_D \otimes D 
% \in \PPDiv_\QQ(Y,\sigma).
$,
then we will always assume that
$\bigcup_{\Delta_D = \emptyset} \supp D$ 
is the support of an effective, semiample
divisor, 
and we understand $\mathcal{D}\in \PPDiv_\QQ(Y,\sigma)$ as
$\mathcal{D}|_{\Loc(\mathcal{D})} \in \PPDiv_\QQ(\Loc(\mathcal{D}),\sigma)$
with 
$$
\Loc(\mathcal{D}):= Y\setminus \bigcup_{\Delta_D = \emptyset} \supp D.
$$
This new convention is compatible with the following evaluation of the 
coefficients of a polyhedral divisor.%
}%\nenh

\begin{definition}
\label{def-wspc}
Let $N$ be a lattice, $\sigma \subseteq N_\QQ$ 
a pointed polyhedral cone,
and $\mathcal{D} = \sum \Delta_D \otimes D$
a polyhedral divisor %with enhanced coefficients
on a normal variety $Y$.
If
$$
\mu \colon
\{\mbox{prime divisors on $Y$}\} \ \to \ \RR
$$
is any map, then
we define the associated 
{\em weighted sum of the polyhedral coefficients\/}
to be
$$
\Delta_\mu 
\ := \ \D_\mu
\ := \
\mu(\mathcal{D})
\ := \
\sum \mu(D) \cdot \Delta_D 
\ \in \ 
\enh{\EnPol_\sigma(N).}
\nenh{\Pol_\sigma(N).}
$$
\end{definition}

% \correct{
% This setting is used because it covers
% many canonical constructions; here are 
% some typical examples.}{}

\begin{example}
\begin{enumerate}
\item
\label{weigh-tail}
For the trivial map $\mu\equiv 0$,
the weighted sum 
$\Delta_0$ % $\wspc_0(\mathcal{D})$ 
gives the common tail cone $\tail(\mathcal{D})$
of the coefficients of $\mathcal{D}$.
\item
\label{weigh-coef}
Fixing a prime divisor $P \in \WDiv(Y)$, 
we may consider $\mu_P(D):=\delta_{D,P}$.
The corresponding 
$\D_P=\mu_P(\mathcal{D})$ % $\wspc_i(\mathcal{D})$ 
recovers the coefficient $\Delta_{P}$ of
$P$.
\item
\label{weigh-local}
Given a point $y\in Y$, 
set $\mu_y(D) := 1$ 
if $y \in D$
and $\mu_y(D) := 0$ else.
Then 
$\mu_y(\mathcal{D})$ 
is precisely the fiber polyhedron
$\Delta_y$ of the point $y  \in Y$.
\item
\label{weigh-deg}
% If $Y$ is a curve and
% $\mu(D)= 1$ for all $D$, then
% $\Delta_{\mu}$ 
% is a polyhedron satisfying
% $\min \bangle{u, \Delta_{\mu}} = \deg \mathcal{D}(u)$
% for all $u\in \tail(\mathcal{D})^\vee$.
% This suggests the notion 
% of a polyhedral degree:
% $\deg \mathcal{D} := \Delta_{\mu}$.
% \item
\label{weigh-intersect}
\correct{}{%
If $C\subseteq Y$ is a curve, then $\mu_C(D):=(C\cdot D)$
% provides a generalization of the previous example.
% The resulting weighted sum may be denoted as 
leads to
$\mu_C(\mathcal{D})=:(C\cdot\mathcal{D})\in \Pol_\sigma(N)$.
In the case of $Y=C$, or
$Y=\PP^n$ and $C$ being the line, we denote
$(C \cdot \D)$ also by $\deg \D$.}
\end{enumerate}
\end{example}

\section{Open embeddings}
\label{sec:opemb1}

In this section, we begin the 
study of open 
embeddings of affine $T$-varieties 
in terms of pp-divisors.
The first statement is a description 
of the equivariant basic open sets obtained by homogeneous
localization.
Recall that in toric geometry equivariant 
localization corresponds to passing 
to a face of a given cone. 
The generalization to pp-divisors involves 
also operations on the base variety; 
here is the precise procedure.
Fix a lattice $N$ and a normal, semiprojective 
variety $Y$. Moreover, let 
$\sigma \subseteq N_\QQ$
be a pointed cone and consider 
a pp-divisor 
$$
\mathcal{D}
\ = \ 
\sum \Delta_D \otimes D
\ \in \ 
\PPDiv_{\QQ}(Y,\sigma).
$$
As usual, $\mathcal{A}$ denotes the 
associated sheaf of $M$-graded 
$\mathcal{O}_Y$-algebras, 
$A := \Gamma(Y, \mathcal{A})$
is the algebra of global sections,
and we set $\t{X} := \Spec_Y(\mathcal{A})$,
and $X := \Spec(A)$.

\begin{definition}
\label{def:loc}
Let $w \in \sigma^\vee \cap M$ and 
$f \in A_w = \Gamma(Y,\mathcal{O}(\mathcal{D}(w)))$.
\begin{enumerate}
\item 
The {\em face\/} of $\Delta \in \Pol_\sigma^+(N)$ 
defined by $w$ is 
\vspace{-1ex}
% \begin{eqnarray*} 
% \face(\Delta,w)
% & := &
% \{v \in \Delta; \; 
% \bangle{w,v} \le  \bangle{w,v'} 
% \text{ for all } v' \in \Delta\}
% \\
% & \in &
% \Pol_{\sigma \cap w^\perp}^+(N).
% \end{eqnarray*}
$$ 
\face(\Delta,w)
 := 
\{v \in \Delta; \; 
\bangle{w,v} \le  \bangle{w,v'} 
\text{ for all } v' \in \Delta\}
\in 
\Pol_{\sigma \cap w^\perp}^+(N).
$$
\item
The {\em zero set\/} of $f$ and the 
{\em principal set\/} 
associated to $f$ are
\vspace{-1ex}
$$ 
Z(f) \ := \ \Supp(\div(f) + \mathcal{D}(w)),
\quad
Y_f \ := \ Y \setminus Z(f).
$$
\item
The {\em localization\/} of the pp-divisor $\mathcal{D}$
by $f$ is
\vspace{-1ex}
% \begin{eqnarray*} 
% \mathcal{D}_f
% &  := &
% \sum \face(\Delta_D, w) \otimes D_{\vert Y_f}.
% \\
% & \in & 
% \CDiv(Y_f,\sigma \cap w^\perp).
% \end{eqnarray*}
$$
\mathcal{D}_f
  := 
\sum \face(\Delta_D, w) \otimes D_{\vert Y_f}
= \emptyset \otimes (\div(f) + \mathcal{D}(w)) 
+ \sum \face(\Delta_D,w) \otimes D.
%  \in 
% \CDiv(Y_f,\sigma \cap w^\perp).
$$
\end{enumerate}
\end{definition}

\begin{lemma}
\label{lastlemma}
Let $w \in \sigma^\vee \cap M$ and 
$f \in A_w = \Gamma(Y,\mathcal{O}(\mathcal{D}(w)))$
as in Definition~\ref{def:loc}.
Then, for $u \in \sigma \cap w^\vee$ and $k \gg 0$,
one has $u+kw \in \sigma^\vee$ and
\begin{eqnarray*}
\mathcal{D}_f(u) 
& = & 
\mathcal{D}(u+kw)|_{Y_f} - \mathcal{D}(kw)|_{Y_f}.
\end{eqnarray*}
\end{lemma}

\begin{proof}
Set $\sigma_w := \sigma \cap w^\perp$
and $\Delta_w := \face(\Delta,w)$.
The first part of the assertion is clear by
$\sigma_w = \sigma^\vee - \QQ_{\ge 0}w$.
The second part is obtained by 
comparing the non-empty coefficients of 
the prime divisors.
For $\mathcal{D}_f(u)$, they are of the form
$\min \langle \Delta_w,u \rangle$.
If $u$ attains this minimum at
$v \in \Delta_w$, then $v$
provides a minimal value for 
$u+kw$ on the whole $\Delta$.
Thus, the claim follows from
$$
\min \langle \Delta_w,u\rangle
\; = \;
\langle v,u \rangle
\; = \;
\langle v,u+kw \rangle - \langle v,kw \rangle
\; = \;
\min \langle \Delta, u+kw \rangle -
\min \langle \Delta, kw \rangle.
$$
\end{proof}

\begin{proposition}
\label{prop:localization}
For a pp-divisor $\mathcal{D}$
on a normal, semiprojective variety $Y$,
let $\mathcal{D}_f$ be the localization
of $\mathcal{D}$ by a homogeneous 
$f \in A_w$.
Then 
% $Y_f$ is semiprojective, 
$\mathcal{D}_f$ is a pp-divisor 
on $Y_f$, 
and the canonical map 
$\mathcal{D}_{f} \to \mathcal{D}$
describes the open embedding 
$X_{f} \to X$.
\end{proposition}

\begin{proof}
% [Proof of Proposition~\ref{prop:localization}]
%
% We first show that $Y_f$ is semiprojective. 
% If $D := \div(f) + \mathcal{D}(w)$ is ample, 
% then this is obvious.
% If $D$ is only semiample, then we have 
% $D = p^*D'$ with an effective ample divisor 
% $D'$ on a semiprojective variety $Y'$ 
% and a projective morphism
% $p \colon Y \to Y'$ having connected fibers.
% Clearly, the semiprojectivity of 
% $Y' \setminus \Supp(D')$ 
% implies semiprojectivity of 
% $Y \setminus \Supp(D)$.
%
% For the remaining claims, 
\nenh{We may assume that $\mathcal{D}$ has non-empty coefficients.}
% We need
% a preparatory consideration.
Recall that $Y_{f}$ is obtained 
by removing the support of
$D = \div(f)+\mathcal{D}(w)$
from $Y$.
In particular, $\mathcal{D}(w)$ is principal on 
$Y_{f}$, and thus, for $k \gg 0$, 
Lemma~\ref{lastlemma} gives
$$
\mathcal{O}_{Y_{f}}(\mathcal{D}_f(u))
\ \cong \
\mathcal{O}_{Y_{f}}(\mathcal{D}(u+kw)).
$$
Using this, one sees that the assignment 
$u \mapsto \mathcal{D}_f(u)$ inherits 
from $u \mapsto \mathcal{D}(u)$ the 
properties~(i) to~(iii) of a pp-divisor 
formulated in Section~\ref{sec:affcase}.
%in~\cite[Def.~2.10]{AlHa} .

To see that $\mathcal{D}_{f} \to \mathcal{D}$
defines an open embedding 
$X_{f} \to X$, 
it suffices to verify that, for any 
linear form 
$u \in (\sigma \cap w^\perp)^\vee\cap M 
= (\sigma^\vee\cap M) - \NN \cdot w$,
we have
$$
\bigcup_{k\gg0}\,\Gamma \big(Y,\,\mathcal{D}(u+kw)\big) / f^k
\; = \;
\Gamma\big(Y_{f}, \,\mathcal{D}(u+kw)- k\, \mathcal{D}(w)\big)
\quad \text{where } k \gg 0.
$$

Consider an element $g/f^k$ of the left hand side. 
Then $\div(g)+\mathcal{D}(u+kw)\ge 0$ holds.
Hence, still on $Y$, we have
$$
\div(g/f^k) +\mathcal{D}(u+kw) - k\,\mathcal{D}(w)
\; \geq \; 
- \div(f^k) - k\,\mathcal{D}(w)
\; = \; -k\,Z(f).
$$
Thus,
$\div(g/f^k) + \mathcal{D}(u+kw) -k\,\mathcal{D}(w)$
is effective on $Y_{f}$, which means that $g/f^{k}$
belongs to the right hand side.

For the reverse inclusion, take any element 
from the right hand side; we may write this 
element as $g/f^k$ with $k\gg 0$.
{From} the relation
$$
\div(g/f^k)+\mathcal{D}(u+kw)-k \,\mathcal{D}(w)
\; \geq \; 
0
$$ 
on $Y_{f}$, we obtain the existence of an $\ell\in\ZZ$ 
such that the same divisor is
$\geq-\ell\,Z(f)$ on $Y$. Moreover, we may assume that 
$\ell\geq k$.
Then,
$$
\div(g/f^k) + \mathcal{D}(u+kw) - k \,\mathcal{D}(w)
\; \geq \;
-\div(f^\ell) -\ell \,\mathcal{D}(w)
$$
holds on $Y$. Using the convexity property of the 
assignment
$u \mapsto \mathcal{D}(u)$, we can conclude
$$
\div(gf^{\ell-k})+\mathcal{D}(u+\ell w) 
\; \ge \;
\div(gf^{\ell-k})+\mathcal{D}(u+kw)+ (\ell-k)\,\mathcal{D}(w)
\;\geq \; 0.
$$
However, this shows that $g/f^k=gf^{\ell-k}/f^\ell$ belongs to the big union
of the left hand side.
\end{proof}

\enh{%
\begin{remark}
\label{rem:enhloc}
Localization can easily be transferred into 
the language of enhanced coefficients.
Let
$$
\mathcal{D}
\ = \ 
\sum \Delta_D \otimes D
\ \in \ 
\EnPPDiv_{\QQ}(Y,\sigma).
$$
For a given homogeneous $f \in A_w$,
let $D(f) := \div(f) + \mathcal{D}(w)$.
Then the localization of $\mathcal{D}$
is
$$
\mathcal{D}_f
\ := \ 
\emptyset \otimes D(f) + \sum \face(\Delta_D,w) \otimes D  
\ \in \ 
\EnPPDiv_{\QQ}(Y,\sigma_w).
$$
\end{remark}
}%\enh

Whereas in toric geometry every equivariant
open embedding of affine toric varieties is 
a localization, this needs no longer hold
for general $T$-varieties.
Thus, in view of equivariant gluing, 
we have to take care of more general
affine open embeddings.
We consider the following situation.
By $N$, we denote again a lattice, 
and $Y$ is a normal variety.
Moreover, 
$\sigma' \subseteq \sigma \subseteq N_\QQ$
are pointed polyhedral cones,
and we consider two pp-divisors
\enh{%
$$
\mathcal{D}' 
\ = \
\sum \Delta'_D \otimes D 
\ \in \
\EnPPDiv_\QQ(Y,\sigma'),
\quad 
\mathcal{D} 
\ = \
\sum \Delta_D \otimes D 
\ \in \
\EnPPDiv_\QQ(Y,\sigma)
$$
with enhanced coefficients.
}%\enh
\nenh{%
$$
\mathcal{D}' 
\ = \
\sum \Delta'_D \otimes D 
\ \in \
\PPDiv_\QQ(Y,\sigma'),
\quad 
\mathcal{D} 
\ = \
\sum \Delta_D \otimes D 
\ \in \
\PPDiv_\QQ(Y,\sigma).
$$
}%\nenh
We suppose that $\Delta'_D \subseteq \Delta_D$
holds for every prime divisor $D \in \WDiv(Y)$. 
For the respective loci of these divisors, we then 
obtain
$$ 
V' 
\ := \
\Loc(\mathcal{D}') 
\ = \ 
\{y \in Y; \; \Delta'_y \ne \emptyset\}
\ \subseteq \
\{y \in Y; \; \Delta_y \ne \emptyset\}
\ = \ 
\Loc(\mathcal{D}) 
\ =: \
V.
$$
Note that we have a natural map 
$\mathcal{D}' \to \mathcal{D}$
of pp-divisors.
As mentioned in Section~\ref{sec:affcase},
this gives rise to a commutative diagram
of $T$-equivariant morphisms,
where the rows contain the geometric 
data associated
to $\mathcal{D}$ and $\mathcal{D}'$
respectively:
$$ 
\xymatrix{
V 
&
{\t{X}} \ar[l]_{\pi} \ar[r]^{r}
&
X
\\
V' \ar[u] 
&
{\t{X}'} \ar[l]^{\pi'} \ar[r]_{r'} \ar[u] 
&
X' \ar[u] .
}
$$
%Moreover, we shall denote by $\mathcal{A}$,
%$\mathcal{A}'$ the sheaves of algebras
%associated to $\mathcal{D}$, $\mathcal{D}'$,
%and by $A$, $A'$ the respective rings of global
%sections.

\begin{proposition}
\label{openemb}
The morphism $X' \to X$ associated to 
$\mathcal{D}' \to \mathcal{D}$ is an open
embedding if and only if any 
$y \in V'$ admits 
$\ku \in \sigma^\vee \cap M$
and 
$f \in A_\ku$
with
$$
y \in V_f \subseteq V',
\quad 
\Delta'_y = \face(\Delta_y,\ku),
\quad
\face(\Delta'_v,\ku) = \face(\Delta_v,\ku)
\text{ for every } v \in V_f.
$$ 
\end{proposition}

\begin{proof}[Proof]
Suppose that $X' \to X$ is an open
embedding.
For short we write $X' \subseteq X$.
Given $y \in V'$, let
$T \mal z' \subseteq (\pi')^{-1}(y)$ 
the (unique) closed $T$-orbit and
choose $f \in A_\ku$,
where $\ku \in \sigma^\vee \cap M$,
such that
$$ 
f(r'(z')) \ne 0,
\qquad
f_{\vert X \setminus X'} = 0.
$$
Then we always have $X_{f} = X'_{f}$.
Since the maps $r \colon \t{X} \to X$ 
and $r' \colon \t{X}' \to X'$ are 
birational and proper, this implies
$$
B 
\ := \ 
\Gamma(\t{X}_{f},\mathcal{O}) 
\ = \
\Gamma(X_{f},\mathcal{O}) 
\ = \
\Gamma(X'_{f},\mathcal{O}) 
\ = \
\Gamma(\t{X}'_{f},\mathcal{O}) 
\ =: \
B'.
$$
Considering the invariant parts, 
we obtain that
$\Gamma(V_{f},\mathcal{O})$
equals 
$\Gamma(V'_{f},\mathcal{O})$.
Since both $V'_f\subseteq V_f$ are semiprojective,
this gives $V'_{f} = V_{f}$. 
Using~\cite[Thm.~3.1~(iii)]{AlHa},
we arrive at  
$$
y 
\ \in \ 
\pi(r^{-1}(X_f)) 
\ = \
V_{f}
\ = \ 
V'_{f}
\ \subseteq \ V'
.
$$
Moreover, $B$ and $B'$ are 
the algebras of global sections 
of the localized pp-divisors 
$\mathcal{D}_f$ and 
$\mathcal{D}'_f$ living
on $V_f = V'_f$.
By~\cite[Lemma.~9.1]{AlHa},
$B = B'$ implies 
$\mathcal{D}_f = \mathcal{D}'_f$.
Thus, we obtain
$$ 
\face(\Delta'_v,\ku) 
\ = \ 
\face(\Delta_v,\ku)
$$
for every $v \in V_f$.
Finally, $f(r'(z')) \ne 0$ implies
$\Delta'_y = \face(\Delta'_y,\ku)$,
which, together with the preceding
obervation, shows
$\Delta'_y = \face(\Delta_y,\ku)$.

Now suppose that 
$\mathcal{D}'$ and
$\mathcal{D}$ 
satisfy the assumptions
of the proposition.
For every $y \in V'$ choose
$\ku$ and $f \in A_\ku$ 
as in the assertion.
Then we have
$\Delta'_y = \face(\Delta'_y,\ku)$.
{From} this,
we can conclude
$(\pi')^{-1}(y) \subseteq \t{X}'_f$.
Consequently, the sets 
$\t{X}'_f$, where  $y \in V'$,
cover $\t{X}'$.

Moreover, the assumption
implies that the  
localized pp-divisors 
$\mathcal{D}'_f$ and 
$\mathcal{D}_f$ 
coincide.
Hence, the canonical maps
$\t{X}'_{f} \to \t{X}_{f}$
are isomorphisms,
and this also holds for the 
canonical maps
$X'_f \to X_f$.
Since $\t{X}'$ is covered 
by the sets 
$\t{X}'_f$, 
we obtain that $X' \to X$ 
is an open embedding.
\end{proof}

\begin{remark}
\label{rem-openface}
Suppose we are in the situation 
of Proposition~\ref{openemb}.
\begin{enumerate}
\item
The condition
$\face(\Delta'_v,\ku) = \face(\Delta_v,\ku)$
for every $v \in V_f$ 
is equivalent to the following one:
% $\face(\Delta'_D,\ku) \ne \face(\Delta_D,\ku)$
% holds at most for prime divisors
% $D \in \WDiv(Y)$ supported in $Z(f)$.
If $\face(\Delta'_D,\ku) \ne \face(\Delta_D,\ku)$, then
$D \in \WDiv(Y)$ is a prime divisor supported in $Z(f)$.
\item
The condition of the previous Proposition~\ref{openemb}
implies that $\Delta'_y \preceq \Delta_y$ holds for all
$y \in Y$. This weaker condition turns out to be equivalent
to the map $\t{X}' \to \t{X}$ 
% associated to $\mathcal{D}' \to \mathcal{D}$
being an open embedding.
\end{enumerate}
\end{remark}

\begin{example}
\label{ex-blowup}
Let $Y=\PP^1$ and $N=\ZZ$. The pp-divisor
$\D=[0,\infty)\otimes \{0\} + [1,\infty)\otimes \{\infty\}$ describes
$\KK^2$ with its standard $\KK^\ast$-action.
On the other hand, we may consider
$\D':= [0,\infty)\otimes \{0\} + \emptyset\otimes \{\infty\}$.
The morphism $\D'\to\D$ describes the blowing up of the origin in $\KK^2$,
hence, it is not an open embedding.
\end{example}

\section{Patchworking}
\label{sec:opemb2}

In this section, we continue the 
study of equivariant open embeddings.
Given a pp-divisor and its associated 
affine $T$-variety $X$, our aim is 
to construct a pp-divisor 
for an invariant, affine, open subset 
$X' \subseteq X$.
Clearly, $X'$ is a union of homogeneous
localizations of $X$. 
We need the following setting.

\begin{definition}
\label{def:reduce}
Let $X$ be an affine $T$-variety,
and let $X' \subseteq X$ be an
invariant, open,  affine subset. 
We say that 
$f_1, \ldots, f_r \in \Gamma(X,\mathcal{O})$
{\em reduce\/} $X$ to $X'$ if
\begin{enumerate}
\item
each $f_i$ is homogeneous and
$X' = \bigcup_{i=1}^r X_{f_i}$
holds,
\item
each $f_i$ is invertible
on some orbit closure in $X'$.
\end{enumerate}
\end{definition}

\begin{remark}
For any invariant, affine, open subset 
 $X' \subseteq X$ of an affine $T$-variety
$X$, there exist homogeneous functions 
$f_1, \ldots, f_r \in \Gamma(X,\mathcal{O})$
that reduce $X$ to $X'$.
If $X' \hookrightarrow X$  is an open embedding that
arises from a map of
pp-divisors $\mathcal{D}' \to \mathcal{D}$
as in Proposition~\ref{openemb},
\correct{}{then the functions $f \in A_\ku$ mentioned there
will do.}
\end{remark}

% \begin{remark}
% If $X' \to X$  is an open embedding that
% arises from a map of
% pp-divisors $\mathcal{D}' \to \mathcal{D}$
% as in Proposition~\ref{openemb},
% then the functions $f \in A_\ku$ mentioned there
% reduce $X$ to (the image of) $X'$.
% \end{remark}

By~\cite[Thm.~8.8]{AlHa},
the open embedding $X'\hookrightarrow X$ may be represented 
by some map of pp-divisors $\mathcal{D}' \to \mathcal{D}$.
In the following, we will show that $\mathcal{D}'$ may be chosen
to live on the same base $Y$ as $\mathcal{D}$ does.

\begin{proposition}
\label{opensubset}
Consider a pp-divisor 
$\mathcal{D} 
= \sum_D \Delta_D \otimes D$
on a normal semiprojective 
variety $Y$, denote the 
associated geometric 
data by
$$\xymatrix{
Y 
& 
{\t{X}} \ar[l]_{\pi} \ar[r]^{r}
&
X
},
$$
and let $X' \subseteq X$ be an invariant,
affine, open subset.
Then $Y' := \pi(r^{-1}(X')) \subseteq Y$ is open 
and semiprojective.
Moreover, if $f_i \in A_{\ku_i}$ reduce $X$ to
$X'$, then
$$
\mathcal{D}' 
\ := \  
\bigcup \mathcal{D}_{f_i}
\ := \
\sum \Delta'_D \otimes D_{\vert Y'},
\quad
\text{where }
\Delta'_D 
\ := \ 
\bigcup_{D \cap Y'_{f_i} \ne \emptyset} 
\face(\Delta_D,\ku_i)
\ \preceq \
\Delta_D
$$
is a pp-divisor on $Y' = \bigcup Y_{f_i}$, 
and
the canonical map
$\mathcal{D}' \to \mathcal{D}$ 
defines an open embedding of
affine varieties having 
$X'$ as its image.
\end{proposition}

\enh{%
Like localization, see~Remark~\ref{rem:enhloc},
the above statement can easily be put into the 
framework of enhanced coefficients.
For later applications, we note in this setting
the following consequence.
}%\enh

\begin{corollary}
\label{sec:cor-intersect-divisors}
Let $X$ be the affine variety arising from
a pp-divisor $\mathcal{D}$ 
\enh{%
with enhanced coefficients 
}%\enh
on a normal variety $Y$, 
and let the pp-divisors
$$ 
\mathcal{D}'
\ = \ 
\sum \Delta'_D \otimes D
\ = \ 
\bigcup \mathcal{D}_{f_i}, 
\qquad
\mathcal{D}''
\ = \ 
\sum \Delta''_D \otimes D
\ = \ 
\bigcup \mathcal{D}_{g_j}
$$
with loci $Y' \subseteq Y$ and $Y'' \subseteq Y$,
respectively,
describe open subsets $X' = \bigcup X_{f_i}$ 
and $X'' = \bigcup X_{g_j}$ as in 
Proposition~\ref{opensubset}.
Then we have
$$ 
\mathcal{D}' \cap \mathcal{D}''
\ := \ 
\sum (\Delta'_D \cap \Delta''_D) \otimes D
\ = \ 
\bigcup \mathcal{D}_{f_ig_j}. 
$$
In particular, 
$\mathcal{D}' \cap \mathcal{D}''$
is a pp-divisor with locus $Y' \cap Y'' \subseteq Y$, 
and the canonical map 
$\mathcal{D}' \cap \mathcal{D}'' \to \mathcal{D}$
describes an open embedding having
$X' \cap X''$ as its image.
\end{corollary}

For the proof of Proposition~\ref{opensubset},
we need two preparatory lemmas.
Let $Y$, $Y''$ be normal 
semiprojective varieties, 
$N$, $N''$ lattices, 
$\sigma \subset N_\QQ$ and
$\sigma'' \subset N''_\QQ$
pointed cones, 
and consider pp-divisors
$$
\mathcal{D} 
\ = \
\sum \Delta_D \otimes D  
\ \in \
\PPDiv_{\QQ}(Y, \sigma),
\quad
\mathcal{D}'' 
\ = \
\sum \Delta''_D \otimes D''  
\ \in \
\PPDiv_{\QQ}(Y'',\sigma'')
$$ 
\nenh{%
with non-empty coefficients.
}%\nenh
Moreover, let 
$(\psi,F,\mathfrak{f})$ be a
map from $\mathcal{D}''$ to
$\mathcal{D}$.
As indicated in 
Section~\ref{sec:affcase},
the map 
$(\psi,F,\mathfrak{f})$ 
gives rise to a commutative diagram
of equivariant morphisms,
where the rows contain the geometric 
data associated
to $\mathcal{D}$ and $\mathcal{D}''$
respectively:
$$ 
\xymatrix{
Y 
&
{\t{X}} \ar[l]_{\pi} \ar[r]^{r}
&
X
\\
Y'' \ar[u]^{\psi}
&
{\t{X}''} \ar[l]^{\pi''} \ar[r]_{r''} \ar[u]^{\t{\varphi}} 
&
X'' \ar[u]_{\varphi} .
}
$$

\begin{lemma}
\label{openemb2}
In the above notation,
suppose that the morphism
$\varphi \colon X'' \to X$ 
is an open embedding.
Then the following holds.
\begin{enumerate}
\item
We have $\t{\varphi}(\t{X}'') = r^{-1}(\varphi(X''))$,
and the induced morphism
$\t{\varphi} \colon \t{X}'' \to \t{\varphi}(\t{X}'')$
is proper and birational.
\item
The image $\psi(Y'') \subseteq Y$ 
is open and semiprojective, and 
$\psi \colon Y'' \to \psi(Y'')$ is 
a projective birational morphism.
\item 
For every $y \in \psi(Y'')$, 
the intersection 
$\t{U}_y := \pi^{-1}(y) \cap \t{\varphi}(\t{X}'')$
contains a unique $T$-orbit that is
closed in $\t{U}_y$.
\end{enumerate}
\end{lemma}

\begin{proof}
Consider the open subset $U := \varphi(X'')$,
its inverse image $\t{U} := r^{-1}(U)$
and $V := \pi(\t{U})$.
By~\cite[Lemma~2.1]{Ha}, the latter set 
is open in $Y$.
These data fit into the 
commutative diagram
$$ 
\xymatrix{
V 
&
{\t{U}} \ar[l]_{\pi} \ar[r]^{r}
&
U
\\
Y'' \ar[u]^{\psi}
&
{\t{X}''} \ar[l]^{\pi''} \ar[r]_{r''} \ar[u]^{\t{\varphi}} 
&
X'' \ar[u]_{\varphi}^{\cong} .
}
$$
The map $\t{\varphi}$ is birational,
because $r''$, $\varphi$ and $r$ are birational.
Moreover, since $r''$ and hence 
$\varphi \circ r''$ are proper,
we infer from the diagram 
that $\t{\varphi}$ is 
proper, and thus surjective.
Consequently, we obtain
$\psi(Y'')=V$;
in particular, this set is 
open in $Y$.
In order to see that $V$ is 
a semiprojective variety, 
note first that for its global
functions, we have
$$ 
\Gamma(V, \mathcal{O})
\ \cong \ 
\Gamma(\t{U}, \mathcal{O})_0
\ \cong \ 
\Gamma(\t{X}'', \mathcal{O})_0
\ \cong \ 
\Gamma(Y'', \mathcal{O});
$$
the first equality is guaranteed 
by~\cite[Lemma~2.1]{Ha}.
Thus, setting 
$A''_0 := \Gamma(Y'', \mathcal{O})$
and $Y''_0 := \Spec(A''_0)$, 
we obtain a commutative diagram
$$ 
\xymatrix{
Y'' \ar[rr]^{\psi} \ar[dr]
& & 
V \ar[dl]
\\
& 
Y''_0
&
}
$$
Since $\psi$ is surjective,
$V \to Y_0''$ is proper.
%(at least for $\KK = \CC$).
Since $V$ is quasiprojective, we 
even obtain that $V \to Y_0''$ is
projective,
and so is 
$\psi \colon Y'' \to V$.
The map $\psi$ is also birational, because 
we have the commutative diagram
$$ 
\xymatrix{
{\KK(Y)}
\ar[r]^{\psi^*}
\ar[d]_{=}
& 
{\KK(Y'')} 
\ar[d]^{=}
\\
{\KK(\t{X})_0}
\ar[r]_{\t{\varphi}^*}^{\cong}
&  
{\KK(\t{X}'')_0}
}
$$
Now, consider $y \in V$ and the intersection
$\t{U}_y := \pi^{-1}(y) \cap \t{\varphi}(\t{X}'')$. 
Since $\pi^{-1}(y)$ contains only finitely 
many $T$-orbits, the same holds for 
$\t{U}_y$.
Let $T \mal z_1, \ldots T \mal z_r$ be the 
closed $T$-orbits of $\t{U}_y$.
We claim 
$$
\psi^{-1}(y)
\ = \
\pi''(\t{\varphi}^{-1}(\t{U}_y))
\ = \
\bigcup_{i=1}^r \pi''(\t{\varphi}^{-1}(T \mal z_i)).
$$
The first equality is clear by surjectivity 
of $\t{\varphi}$ and the quotient maps
$\pi, \pi''$.
The second one is verified below; 
it uses properness of $\t{\varphi}$:
Given $y'' \in \pi''(\t{\varphi}^{-1}(\t{U}_y))$,
we have $y'' = \pi''(z'')$ for some $z'' \in 
\t{\varphi}^{-1}(\t{U}_y)$.
Since $\pi''$ is constant on orbit closures,
we may assume that $T'' \mal z''$ is closed 
in $\t{\varphi}^{-1}(\t{U}_y)$.
By properness of $\t{\varphi}$, the image
$\t{\varphi}(T'' \mal z'') = T \mal \t{\varphi}(z'')$
is closed in $\t{U}_y$. 
It follows that $y''$ belongs to the right hand
side.  

Having verified the claim, we may proceed as follows.
The closed invariant subsets 
$\t{\varphi}^{-1}(T \mal z_i) \subseteq \t{X}''$
are pairwise disjoint. 
By the properties of the good quotient 
$\pi'' \colon \t{X}'' \to Y''$,
the images $\pi''(\t{\varphi}^{-1}(T \mal z_i))$
are pairwise disjoint as well.
In particular, $\psi^{-1}(y)$ 
is disconnected if $r>1$. 
The latter is impossible because $\psi$, 
as a birational projective
morphism between normal varieties, has 
connected fibers.
\end{proof}

\begin{lemma}
\label{faceunion}
For the functions $f_i \in A_{\ku_i}$ 
of Proposition~\ref{opensubset},
we always have 
$$
\Delta'_D 
\ := \ 
\bigcup_{D \cap Y'_{f_i} \ne \emptyset} 
\face(\Delta_D,\ku_i)
\ \preceq \
\Delta_D.
$$
In particular, for every $f_i$ with 
$D \cap Y'_{f_i} \ne \emptyset$,
we have 
$\face(\Delta_D,\ku_i) \preceq \Delta'_D$.
\end{lemma}

\begin{proof}
Let $D$ be a prime divisor intersecting $Y'$, 
and consider a point $y \in D \cap Y'$ 
such that $y \in Y_{f_i}$ holds for all
the $f_i$ of Proposition~\ref{opensubset} 
with $D \cap Y_{f_i} \ne \emptyset$ and
$D$ is the only prime divisor with 
$\Delta_D \ne 0$ containing $y$. 
Then we have 
\begin{eqnarray*}
\Delta_y
& = & 
\Delta_D.
\end{eqnarray*}

According to \cite[Thm.~8.8]{AlHa},
the inclusion $X' \subseteq X$ is
described by a map of pp-divisors.
Thus, we may apply Lemma~\ref{openemb2}
and obtain that there is a unique
closed $T$-orbit $T \mal z$ in 
$\pi^{-1}(y) \cap r^{-1}(X')$.
This orbit corresponds to a face
$F(z) \preceq \Delta_y$ via
$$ 
T \mal z 
\ \mapsto \
\omega(z)
\ \mapsto \
\face(\Delta_y,u)
\text{ with }
u \in \interior \omega(z),
$$
where $\interior \omega(z)$ denotes the relative interior
of the cone $\omega(z)$.
Since $z \in r^{-1}(X')$ holds, 
some of the $f_i$ of Proposition~\ref{opensubset}
satisfies $f_i(z) \ne 0$
but vanishes on 
$\b{T \mal z} \setminus T \mal z$,
where the closure is taken in 
$\pi^{-1}(y)$.
This means $\ku_i \in \interior\omega(z)$.
To conclude the proof, it
suffices to show 
\begin{eqnarray*}
\bigcup_{D \cap Y'_{f_j} \ne \emptyset} 
\face(\Delta_y,\ku_j)
& = & 
\face(\Delta_y,\ku_i).
\end{eqnarray*}

For any $f_j$ with $D \cap Y'_{f_j} \ne \emptyset$,
we have $y \in Y_{f_j}$.
Thus, there is a point $z_j \in \pi^{-1}(y) \cap X'$
with $f(z_j) \ne 0$.
We may choose $z_j$ such that $\omega(z_j)$
is minimal; this means $\ku_j \in \interior\omega(z_j)$.
Since $z \in \b{T \mal z_j}$ holds,
we obtain $\omega(z) \preceq \omega(z_j)$.
This in turn implies
$\face(\Delta_y,\ku_j) \preceq \face(\Delta_y,\ku_i)$,
and the above equation follows.
\end{proof}

\begin{proof}[Proof of Proposition~\ref{opensubset}]
By~\cite[Thm.~8.8]{AlHa},
there is a map $(\psi,F,\mathfrak{f})$ 
of pp-divisors 
$\mathcal{D}'' \to \mathcal{D}$
such that the associated morphism 
% $X'' \to X$ equals 
is the open embedding of $X'$ into $X$.
Lemma~\ref{openemb2} gives
$Y' = \psi(Y'')$, and it ensures
openness and semiprojectivity of 
$Y' \subseteq Y$.
Lemma~\ref{faceunion} tells us
$$
\Delta'_D 
\ := \ 
\bigcup_{i=1}^r \face(\Delta_D,\ku_i)
\ \preceq \
\Delta_D.
$$
Property~\ref{def:reduce}~(ii) ensures
that any $X_{f_i}$ contains 
the generic orbit closure of $X'$.
Consequently, all $X_{f_i}$ have the 
same weight cone.
Hence, the tail cone of 
$\face(\Delta_D,\ku_i)$ 
does not depend on $i$,
and, thus, $\mathcal{D}'$ 
is a well defined polyhedral 
divisor on $Y'$.

In order to verify the pp-properties
for $\mathcal{D}'$, we only have 
to concern ourselves with semiampleness and 
bigotry. 
Locally,  we have canonical 
isomorphisms
$$ 
\psi^* (\mathcal{D}')
\ = \ 
\psi^* (\mathcal{D}'_{f_i})
\ = \ 
\psi^* (\mathcal{D}_{f_i})
\ = \ 
\mathcal{D}''_{\psi^* f_i}
\ = \ 
\mathcal{D}''.
$$
This shows $\psi^* (\mathcal{D'}) \cong \mathcal{D}''$,
and thus we obtain both properties by pushing
forward suitable global sections.
%(Hier waere ein direktes Argument natuerlich viel besser.)

Finally, the fact that the map 
$\mathcal{D}' \to \mathcal{D}$ 
is an open embedding with image
$X'$ follows by comparing the 
induced maps
$\mathcal{D}'_{f_i} \to \mathcal{D}_{f_i}$
of the localizations.
\end{proof}

\section{Divisorial Fans}
\label{glue}

In toric geometry, the equivariant 
gluing of affine 
pieces is described by means of a fan,
i.e.\ a collection of polyhedral cones
satisfying natural compatibility conditions.
In this section, we generalize this 
idea and present a natural concept
to describe the equivariant gluing 
of affine varieties with torus action.

\begin{definition}
\label{def:ppface}
Let $N$ be a lattice, $\sigma', \sigma \subseteq N_\QQ$ 
pointed cones, $Y$ a normal, semiprojective variety,
and consider two pp-divisors 
\enh{%
with enhanced coefficients 
}%\enh
on $Y$:
\nenh{%
$$
\mathcal{D}'
 =  
\sum \Delta'_D \otimes D 
\ \in \
\PPDiv_\QQ(Y,\sigma'),
\quad 
\mathcal{D} 
 = 
\sum \Delta_D \otimes D
\ \in \
\PPDiv_\QQ(Y,\sigma).
$$
}%\nenh
\enh{%
$$
\mathcal{D}'
 =  
\sum \Delta'_D \otimes D 
\ \in \
\EnPPDiv_\QQ(Y,\sigma'),
\quad 
\mathcal{D} 
 = 
\sum \Delta_D \otimes D
\ \in \
\EnPPDiv_\QQ(Y,\sigma).
$$
}%\enh
We call $\mathcal{D'}$ a {\em face\/} 
of $\mathcal{D}$
(written $\mathcal{D}' \preceq \mathcal{D}$)
if $\Delta'_D \subseteq \Delta_D$ holds for all $D$
%holds for every prime divisor $D \in \WDiv(Y)$, 
and for any $y \in \Loc(\mathcal{D}')$ there 
are $\ku_y \in \sigma^\vee \cap M$
and a $D_y$ in the linear system $|\D(\ku_y)|$
with
\begin{enumerate}
\item 
$y\notin \supp (D_y)$,
\item
$\Delta'_y = \face(\Delta_y,\ku_y)$,
\item
$\face(\Delta'_v,\ku_y) = \face(\Delta_v,\ku_y)$
for every $v \in Y \setminus \supp(D_y)$.
\end{enumerate} 
\end{definition}

% \begin{remark}
% \begin{enumerate}
% \item
% If a pp-divisor $\mathcal{D}'$ is a face of a 
% pp-divisor $\mathcal{D}$, then the tail cone 
% of  $\mathcal{D}'$ is a face of that of 
% $\mathcal{D}$.
% \item
% Let $\mathcal{D}', \mathcal{D}$
% be pp-divisors on $Y$ such that
% $\Delta'_D \subseteq \Delta_D$ for
% all their coefficients.
% Then Proposition~\ref{openemb} says
% that
% $\mathcal{D}'$ is a face of 
% $\mathcal{D}$ if and
% only if the map $\mathcal{D}' \to \mathcal{D}$
% defines an open embedding
% $X(\mathcal{D}') \to X(\mathcal{D})$.
% \end{enumerate}
% \end{remark}
% Having defined a face relation on pp-divisors,
% we are ready\
%
\correct{}{By Proposition~\ref{openemb},
the face relation 
``$\mathcal{D}' \preceq \mathcal{D}$''
is the combinatorial
counterpart of the open embeddings after applying 
the functor $X(\kbb)$. It implies that
the tail cone
of  $\mathcal{D}'$ is a face of the tail cone of
$\mathcal{D}$.
This enables us}
to generalize the concept of a fan 
to the setting of pp-divisors.
Recall that in the preceding section,
we introduced the intersection of
two polyhedral divisors 
$\mathcal{D}'  =  
\sum \Delta'_D \otimes D$
and 
$\mathcal{D} =  
\sum \Delta_D \otimes D$
\enh{%
with enhanced coefficients 
}%\enh
on a common variety $Y$ as
\begin{eqnarray*}
\mathcal{D}' \cap \mathcal{D}
& := &
\sum (\Delta'_D \cap \Delta_D) \otimes D.
\end{eqnarray*}

\begin{definition}
\label{def:ppfan}
Let $N$ be a lattice, and $Y$ a normal,
semiprojective variety.
A {\em divisorial fan \/} 
on $(Y,N)$ is a set $\mathcal{S}$ of pp-divisors
\enh{%
$\mathcal{D} \in \EnPPDiv_\QQ(Y,\sigma_{\mathcal{D}})$
}%\enh
\nenh{%
$\mathcal{D} \in \PPDiv_\QQ(Y,\sigma_{\mathcal{D}})$
}%\nenh
with tail cones $\sigma_{\mathcal{D}} \subseteq N_\QQ$
such that
for any two $\mathcal{D}',\mathcal{D} \in \mathcal{S}$
the intersection
$\mathcal{D}' \cap \mathcal{D}$
is a face of both
$\mathcal{D}'$ and $\mathcal{D}$
and, moreover, belongs to $\mathcal{S}$.
\end{definition}

Given a divisorial fan $\mathcal{S} =\{\mathcal{D}^i; \; i \in I\}$
on a normal variety $Y$,
we have the affine
$T$-varieties $X_i := X(\mathcal{D}^i)$,
and for any two $i,j \in I$, 
the $T$-equivariant open embeddings
$$ 
\xymatrix{
X_i
& &
X(\mathcal{D}^i \cap \mathcal{D}^j)
\ar[ll]_{\eta_{ij}}
\ar[rr]^{\eta_{ji}}
& &
X_j.
}
$$
We denote the associated images by  
$X_{ij} := \eta_{ij}(X(\mathcal{D}^i \cap \mathcal{D}^j)) \subseteq X_i$.
Then we have $T$-equivariant isomorphisms
$
\varphi_{ij} := \eta_{ji} \circ \eta_{ij}^{-1}
$
from 
$X_{ij}$ onto $X_{ji}$.

\begin{theorem}
\label{sec:thm-divisorial-fans1}
The affine $T$-varieties $X_i$ and the
isomorphisms 
$\varphi_{ij} \colon X_{ij} \to X_{ji}$
are gluing data.
The resulting space
$$ 
X
\ := \ 
X(\mathcal{S})
\ := \ 
\bigsqcup X_i \big/ \sim,
\quad \text{where }
X_{ij} \ni x \ \sim \ \varphi_{ij}(x) \in X_{ji}, 
$$
is a prevariety with affine diagonal $X \to X \times X$, 
and it comes with a (unique) 
$T$-action such that all canonical maps
$X_i \to X$ are equivariant.
\end{theorem}

\begin{proof}
The only thing to check is that 
the maps $\varphi_{ij} \colon X_{ij} \to X_{ji}$
do in fact define gluing data.
Concretely, this means to verify two things,
namely
$$ 
\varphi_{ij}(X_{ij} \cap X_{ik})
\ = \ 
X_{ji} \cap X_{jk}
% \qquad
\hspace{0.7em}\mbox{and}\hspace{0.7em}
\varphi_{ik} 
\ = \ 
\varphi_{jk} \circ \varphi_{ij}.
$$

The first one of these identities
can be directly deduced from the
following observation: for any triple 
$i,j,k \in I$, there is 
a commutative diagram
$$ 
\xymatrix{
X(\mathcal{D}^i \cap \mathcal{D}^j)
\ar[r]^{\eta_{ij}}
&
X(\mathcal{D}^i)
&
X(\mathcal{D}^i \cap \mathcal{D}^k)
\ar[l]_{\eta_{ik}}
\\
&
X(\mathcal{D}^i \cap \mathcal{D}^j \cap \mathcal{D}^k)
\ar[ul]
\ar[u]
\ar[ur]
& 
}
$$
and Corollary~\ref{sec:cor-intersect-divisors}
yields that $X(\mathcal{D}^i \cap \mathcal{D}^j \cap \mathcal{D}^k)$
is mapped onto the intersection $X_{ij} \cap X_{ik}$
of $X_{ij} = \eta_{ij}(X(\mathcal{D}^i \cap \mathcal{D}^j))$
and  $X_{ik} = \eta_{ik}(X(\mathcal{D}^i \cap \mathcal{D}^k))$.

The second identity then may be verified as an identity
of rational maps: all $X_i$ have the 
same function field, and the pull back maps $\varphi_{ij}^*$
are the identity. 
\end{proof}

\begin{example}
\label{ex-omegaP2}
The figure in the introduction on Page \pageref{intro-figure}
shows a divisorial fan on the projective line 
$\PP^1$ generated by six maximal pp-divisors 
$\mathcal{D}^1, \ldots, \mathcal{D}^6$,
all of them of the form
\begin{eqnarray*}
\mathcal{D}^i
& = & 
\Delta^i_0 \otimes \{0\}
\ + \
\Delta^i_1 \otimes \{1\}
\ + \
\Delta^i_\infty \otimes \{\infty\}.
\end{eqnarray*}
In order to indicate in the figure that a polyhedron 
is a coefficient of the divisor 
$\mathcal{D}^i$, we put the label
``$\mathcal{D}^i$'' on it.
% 
% \bigskip
% 
% 
% \begin{center}
% \input{tangp2.pstex_t}
% \end{center}
%
%
%\vspace{-3ex}
%\begin{figure}[hbtp]
%\centering
%\subfigure[$\{0\}$]{\cotangfannull}  
%\subfigure[$\{\infty\}$]{\cotangfaninfty} 
%\subfigure[$\{1\}$]{\cotangfanone}
% \caption{A fan of polyhedral divisors}
%\label{fig:ex-cot-fan}
%\end{figure}
%\vspace{-3ex}\\
%(To understand the meaning of the labeling, see the comment after
%Proposition \ref{prop-coefan}.)
%
As we will see in Section~\ref{sec:txt-bundles},
the corresponding $T$-variety is the projectivized
cotangent bundle $\PP(\Omega_{\PP^2})$ 
over the projective plane.
\end{example}

We conclude this section with a recipe for producing 
lots of examples using toric geometry.
Firstly we recall from~\cite{AlHa} a toric construction
of pp-divisors.

Let $N'$ be a lattice, and denote by 
$M' := \Hom(N',\ZZ)$ the dual lattice.
Let $\delta \subseteq N'_\QQ$ 
be a pointed polyhedral cone. 
We consider the associated affine toric
variety and its big torus
$$
X'
%\TT\VV(\delta)
\ := \ 
\toric{\delta}:=\Spec \, \KK[\delta^\vee \cap M'],
\qquad
T' 
\ := \ 
\Spec \, \KK[M'].
$$
Suppose that $T \to T'$ is a 
monomorphism of tori arising from
a surjection $\deg \colon M'\to M$
of the respective (character) 
lattices; in the case of 
$N' = \ZZ^n$ and $\delta=\QQ^n_{\ge 0}$, 
this map just fixes multidegrees $\deg(z_i) \in M$ 
for the coordinates
$z_1,\ldots,z_n$ of $X' = \KK^n$.

The aim is to construct a pp-divisor
for the induced $T$-action on $X'$.
We will work in terms of the following
exact sequences:
$$ 
\xymatrix@C=3em@R=3ex{
0 \ar[r] &
M'' \ar[r]^-{p^\ast}&
M' \ar[r]^-{\deg} \ar@/^1pc/[l]_-{t^\ast}&
M \ar[r] \ar@/^1pc/[l]_-{s}&
0\\
0 &
\ar[l] N'' \ar@/_1pc/[r]^-{t}&
\ar[l]_-{p} N' \ar@/_1pc/[r]^-{s^\ast}&
\ar[l]_-{\deg^\ast} N &
0 \ar[l]
}
$$
Here we have, additionally, chosen a section $s$ of
$\deg \colon M' \to M$; this corresponds to a section
$t$ of $p \colon N' \to N''$ via 
$t^\ast = \id_{M'} - s \circ \deg$.

Let $\Sigma''$ be the fan in $N''_\QQ$
that is obtained as the coarsest 
common subdivision of the
images of the faces of $\delta$ under $p$; 
its support
is $|\Sigma''|=p(\delta)$. 
Similarly, we obtain a subdivision
of $\omega := \deg(\delta^\vee)$ inside 
$M_\QQ$.
Then the {\em positive fiber\/}
of $u\in \omega\cap M$ is
$$
\Delta(u)
\ := \ 
\big(\deg^{-1}(u)\cap \delta^\vee\big) - s(u)
\ = \
t^\ast\big( \deg^{-1}(u)\cap \delta^\vee\big) 
\ \subseteq \ 
M''_\QQ.
$$
The normal fans $\Lambda(\Delta(u))$ of these polytopes vary
with the chamber structure of $\omega$. 
Their coarsest common refinement
is exactly the fan $\Sigma''$.

The toric variety $Y''=\toric{\Sigma''}$ 
associated to $\Sigma''$ 
is the Chow quotient of $X'$ by the action of 
$T$. 
The polytopes $\Delta(u)$ correspond to semiample 
divisors on $Y''$, and they are precisely the 
evaluations of the pp-divisor on $Y'$ describing
$X'$ as a $T$-variety, namely
\begin{eqnarray*}
\mathcal{D}'
& := &
\sum_{\varrho \in (\Sigma'')^{(1)}} \Delta_{\varrho} \otimes D_{\varrho},
\end{eqnarray*}
where, by abuse of notation, 
the ray $\varrho\in (\Sigma'')^{(1)}$ is identified with 
its primitive lattice vector $\varrho\in N''$, 
where $D_\varrho$ denotes the 
invariant prime divisor corresponding to
$\varrho$, and, dualizing the previous formula for $\Delta(u)$, 
$$\Delta_{\varrho}:= \Delta_{\varrho}(\delta) :=
\big(p^{-1}(\varrho)\cap \delta\big) - t(\varrho)
\ = \
s^\ast\big( p^{-1}(\varrho)\cap \delta\big) 
\ \subseteq \ 
N_\QQ.
$$

Now, more generally, let $X'$ be a semiprojective,
\correct{}{%
not necessarily affine,%
}  
toric variety arising from a fan $\Sigma'$ in $N'$. 
Then, similarly to the above setup, we may consider the 
coarsest common refinement $\Sigma''$ of all 
images $p(\delta)$ where $\delta \in \Sigma'$.
This is again a fan in $N''_\QQ$, and 
we denote by $Y''$ the associated toric variety.

In addition, we have, for every cone $\delta \in \Sigma'$,
the previous construction $\Sigma''(\delta)$ refining 
all projected faces of $\delta$.
Note that each $\vert \Sigma''(\delta) \vert$ is a union of 
cones of $\Sigma''$ and we may define a polyhedral 
divisor 
\enh{%
with enhanced coefficients 
}%\enh
on $Y''$ by
\begin{eqnarray*}
\mathcal{D}'(\delta)
& := & 
\sum_{\varrho \in (\Sigma'')^{(1)}} \Delta_{\varrho} (\delta)
\otimes D_{\varrho}.
%\qquad \text{where }
%\Delta_{\varrho}
%\ := \ 
%\left\{
%\begin{array}{ll}
%\Delta_{v_\varrho} 
%& \text{if } 
%\varrho \subset \vert (\Sigma'')^{(1)} \vert
%\\
%\emptyset 
%& \text{else. }
%\end{array}
%\right.
\end{eqnarray*}
This pp-divisor is the pullback of the one previously
associated to the affine chart $X(\delta) \subseteq X'$. 
Its locus equals $\toric{\Sigma''\cap|\Sigma''(\delta)|}$
which is a modification of $\toric{\Sigma''(\delta)}$.
Elements 
$\varrho\in(\Sigma'')^{(1)}\setminus |\Sigma''(\delta)|$
lead to $ \Delta_{\varrho}(\delta)=\emptyset$ in a natural way.
The pp-divisors $\mathcal{D}'(\delta)$,
where $\delta \in \Sigma'$, obviously 
fit together to a divisorial fan $\mathcal{S}'$ on $Y'$, 
and this divisorial fan describes the $T$-action on the
toric variety $X'$.

\begin{proposition}
\label{toricppfan}
Let $X \subseteq X'$ be a closed $T$-invariant subvariety
with $X \cap T' \ne \emptyset$,
and let $\imath \colon Y \to Y''$ be the normalization of
the closure of the image of
$X \cap T'$ in $Y''$. 
Then the $\mathcal{D}(\delta) := \imath^*(\mathcal{D}'(\delta))$
fit together to a divisorial fan 
$\mathcal{S}:=\imath^\ast\mathcal{S}'$ on $Y$, 
and this divisorial fan describes the $T$-variety $X$. 
\end{proposition}

We leave the proof of this observation to the reader.
Note that if $X \cap T'$
is given by $T$-homogeneous equations 
$f_i \in \Gamma(T', \mathcal{O})$,
where $i \in I$,
then, multiplying with $\chi^{-s(\deg f_i)}$
shifts them into $\KK[M'']$, and
we obtain $Y$ as the normalized closure in $Y''$
of 
\begin{eqnarray*}
V(\chi^{-s(\deg f_i)} f_i; \;  i\in I)
& \subseteq & 
T''=\Spec\,\KK[M''].
\end{eqnarray*}

%\begin{proposition}
%Let $X$ be variety such that any two
%points $x,x' \in X$ admit a common  
%affine neighbourhood $X \setminus \Supp(D)$
%with an effective Cartier divisor $D$ on $X$.
%Suppose moreover, that an algebraic torus
%$T$ acts effectively on $X$.
%Then the $T$-variety $X$ arises from a divisorial fan. 
%\end{proposition}
% 
%\begin{proof}
%According to~\cite{Ha0}, there is a $T$-equivariant
%closed embedding $X \to X'$ into a smooth toric variety 
%$X'$, where $T$ acts as a subtorus of the big torus
%$T' \subseteq X'$.
%Replacing $X'$, if necessary, with a suitable orbit 
%$T'$-orbit closure, one may even assume 
%$X \cap T' \ne \emptyset$. 
%Thus, the assertion follows from 
%Proposition~\ref{toricppfan}.
%\end{proof}

The previous proposition presents a kind of an algorithm
for constructing a divisorial fan for $T$-varieties that are
equivariantly embedded in a toric variety, e.g.\ in a projective space.
Besides this, we have the following.

\begin{theorem}
\label{sec:thm-divisorial-fans2}
Up to equivariant isomorphism,
every normal variety with an
effective algebraic torus action
arises from a divisorial fan.
\end{theorem}

\begin{proof}
Let $X$ be a normal variety with 
an effective action of an algebraic
torus $T$.
Then Sumihiro's Theorem ensures that~$X$ 
is covered by $T$-invariant open affine 
subvarieties $X_i \subseteq X$.
By~\cite[Theorem~3.4]{AlHa},
each $X_i$ arises from a pp-divisor
$\mathcal{D}^i$ living on 
semiprojective varieties $Y_i$.
Choose projective closures $Y_i \subseteq Y_i'$ 
such that each complement $Y_i' \setminus Y_i$ 
is the support of a semiample divisor.
\enh{%
Then each $\mathcal{D}^i$ 
extends canonically 
to a pp-divisor with enhanced coefficients 
}%\enh
\nenh{%
Then by introducing empty coefficients, each $\mathcal{D}^i$
extends canonically
to a pp-divisor 
}%\nenh
on $Y_i'$ having $Y_i$ as its locus.
Note that $X_i \supseteq X_{i} \cap X_j \subseteq X_i$ induces 
rational maps $Y_{i}' \leftrightarrow Y_j'$.
Via resolving indeterminacies~\cite[Example~II.7.17.3]{hartshorne:ag}, 
we obtain a projective variety $Y$ which
dominates all $Y_i'$ and is compatible 
with these rational maps.
Then, using Proposition~\ref{opensubset}
one sees that the pull back divisors
of the $\mathcal{D}^i$ fit together 
to the desired divisorial fan 
on~$Y$.
\end{proof}

%\begin{remark}
%The previous proposition presents a kind of an algorithm
%of how to construct a divisorial fan for $T$-varieties that are
%equivariantly embedded in a toric variety, e.g.\ in a projective space.\\
%
%Besides this, it is true that {\em every} (separated) variety $X$ with
%an effective torus action arises from a divisorial fan.
%Indeed, after choosing an open, affine,
%$T$-invariant covering $\{U_i\}$, each of the pieces $U_i$
%can be represented by a pp-divisor on some $Y_i$. Using
%Proposition \ref{opensubset} and eliminating the points
%of indeterminacy of the rational maps among the $Y_i$,
%cf.\ \cite[Ex II.7.17.3]{hartshorne:ag}, all participating pp-divisors
%can be pulled back to a common base $(Y,N)$ --
%eventually yielding a divisorial fan.
%\end{remark}

%%%%%%%%%%%%%%%%%%%%%%%%
%%%%%%%%%%%%%%%%%%%%%%%%
%%%%%%%%
%%%%%%%%  Coherent fans
%%%%%%%%
%%%%%%%%%%%%%%%%%%%%%%%%%%
%%%%%%%%%%%%%%%%%%%%%%%%
\section{Coherent fans}
\label{sec:coherent}

In Definition \ref{def:ppface}, we have used our characterization
of open embeddings among affine $T$-varieties
to establish the face relation between pp-divisors. This led to the notion
of divisorial fans in Definition \ref{def:ppfan}.
While these fans describe general $T$-varieties, 
a restriction of this generality will lead to some
simplification. 
To be precise, we will overcome the
% At the moment we still face the 
following problems:
Firstly, in the present definition of a divisorial fan,
one has to check pp-ness not
only for the ``maximal'' polyhedral
divisors, but also for all of their respective
intersections. 
Secondly, our characterization of open embeddings is precise,
but tough to deal with. 

% Eventually, one is asked to imagine 
% a divisorial fan rather as a ``fansy divisor'',
% i.e.\ as a divisor having
% polyhedral subdivisions as coefficients
% (we .
% However, whether given polyhedral complexes lead to a valid
% divisorial fan  
% as the ``fan'' of its coefficients. However, whether 
% those polyhedral subdivisions of $N_\QQ$ lead to a 
% valid fan of pp-divisors does
% depend on the position of the divisors on $Y$.

% These coherent divisorial fans 
% turn out to be especially nice in the
% case of $\dim Y=1$, i.e.\ in the case of 1-codimensional torus actions.
% However, we will first 
% make precise the notion of the coefficient 
% of a fan of pp-divisors.

\begin{definition}
\label{def-wspcfan}
Let $\mathcal{S}=\{\mathcal{D}^i\}$ be a set of polyhedral divisors
on $(Y,N)$.
To any map
$$
\mu \colon \{\mbox{prime divisors on $Y$}\} \to \RR ,
$$
% as in Definition \ref{def-wspc},
we define the 
$\mu$-{\em slice}
$\mathcal{S}_\mu:=\mu(\mathcal{S}):=\{\mu(\mathcal{D}^i)\}$
as the 
\correct{}{%
set of the weighted sums of the $\mathcal{D}^i$-coefficients.
If $\mathcal{S}$ is a divisorial fan, then 
$\mathcal{S}_\mu:=\mu(\mathcal{S})$ even
forms a directed system.%
}
According to the examples following Definition \ref{def-wspc},
we obtain as special cases
$\tail(\mathcal{S}):=\mathcal{S}_0$ and
$\mathcal{S}_P$, $\mathcal{S}_y$ for prime divisors $P\subseteq Y$
or points $y\in Y$, respectively.
\end{definition}

% \begin{definition}
% \label{def-wspcfan}
% Let $\mathcal{S}=\{\mathcal{D}^i\}$ be a divisorial fan 
% on $(Y,N)$.
% To any map
% $$
% \mu \colon \{\mbox{prime divisors on $Y$}\} \to \RR ,
% $$
% % as in Definition \ref{def-wspc},
% we define the 
% $\mu$-{\em slice}
% $\mathcal{S}_\mu:=\mu(\mathcal{S}):=\{\mu(\mathcal{D}^i)\}$
% as the directed system
% of the weighted sums of the $\mathcal{D}^i$-coefficients.
% %
% According to the examples following Definition \ref{def-wspc},
% we obtain as special cases 
% $\tail(\mathcal{S}):=\mathcal{S}_0$ and
% $\mathcal{S}_P$, $\mathcal{S}_y$ for prime divisors $P\subseteq Y$
% or points $y\in Y$, respectively.
% \end{definition}

\begin{remark}
\label{rem-coefan}
Let $\mathcal{S}$ be a divisorial fan on $(Y,N)$.
Then, by Remark \ref{rem-openface}(ii),
for any $y\in Y$,
% prime divisor $P$ on $Y$, 
the systems $\CS_y$ form a polyhedral complex,
i.e.\ the transition maps are face relations.
Among them are the slices $\CS_P$ for prime divisor $P$ on $Y$.
They even provide polyhedral subdivisions,
i.e.\ $(\D^1\cap\D^2)_P=\D^1_P\cap\D^2_P$, where the cells are
labeled by the elements of $\CS$. Note that multiple labels
may occur. 
If one of the pp-divisors does not occur as a label in, say,
$\CS_P$,
then this indicates that its polyhedral
$P$-coefficient is empty.
Thus, a divisorial fan can be seen as a ``fansy divisor'' with
the subdivisions $\CS_P$ posing as its coefficients.
Similarily, $\tail(\mathcal{S})$ will be called
the {\em tail fan}.
\end{remark}

\begin{example}
\label{ex-noncoh}
We consider $N=\ZZ$; let $\CS=\{\D^1, \D^2, \D^1\cap\D^2\}$ 
be the divisorial fan
generated by 
$\D^i= \Delta^i_1\otimes D_1 + \Delta^i_2\otimes D_2$ 
with some prime divisors $D_1$, $D_2$ of $Y$ and coefficients
\begin{center}
\input{noncoh.pstex_t}
\end{center}
While the sets $\CS_{D_\nu}$ form nice polyhedral complexes,
the map $\mu$ with $\mu(D_1)=\mu(D_2)=1$ leads to a slice 
$\mu(\CS)$ 
where both $\mu(\D^i)$ coincide with the interval
$[-1,1]$. However, $\D^1\cap\D^2$ has the coefficients
$\Delta^{12}_1=\Delta^{12}_2=\{0\}$, hence 
$\Delta^{12}_\mu=\{0\}$, and this is not a face of $[-1,1]$.
Since $\mu(\D^i)$ is not a complex, it cannot occur as
a slice of the form $\CS_y$. Thus, $D_1\cap D_2=\emptyset$ is a
necessary condition
for $\CS$ to form a divisorial fan.
\end{example}

The next example shows that the slices $\CS_y$,
in contrast to the special case $\CS_P$, need not
 be polyhedral subdivisions.

\begin{example}
\label{ex-noncomp}
Let $Y=\AA^2_\KK$ and denote by $y$ the origin.
Then, $\Delta^1_{y}$ and $\Delta^2_{y}$ are the two
polyhedra in the rightmost figure, but
$\Delta^1_{y}\cap\Delta^2_{y}$ is not a face of the
$\Delta^i_{y}$. However, $\CS_y$ is still a complex since
$\Delta^{12}_y=\emptyset$.

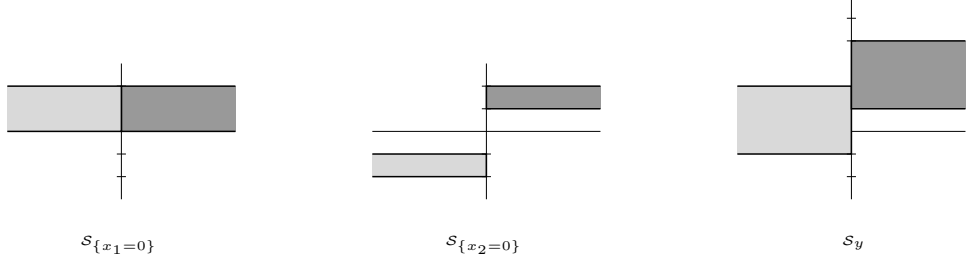
\begin{center}
\input{ex64.pstex_t}
\end{center}

\end{example}

\begin{definition}
\label{def-coherent}
A set of pp-divisors
$\mathcal{S}=\{\mathcal{D}^i=\sum_D \Delta_D^i\otimes D\}$ 
will be called {\em coherent}
if for any $i,j$ there is a $u^{ij}\in M$ 
such that for all $D$ there is a $c_D^{ij}$ with
\vspace{-0.5ex}
$$
\max\langle \Delta^i_D, u^{ij}\rangle 
\leq c_D^{ij}\leq
\min\langle \Delta^j_D, u^{ij}\rangle
\vspace{-0.5ex}
$$
and 
\vspace{-0.5ex}
$$
\Delta^i_D\cap [\langle \kbb, u^{ij}\rangle=c_D^{ij}]=
\Delta^j_D\cap [\langle \kbb, u^{ij}\rangle=c_D^{ij}].
$$
\end{definition}

In principle, coherence means that the coefficients of any two
polyhedral divisors $\D^i$ and $\D^j$ are separated by
hyperplanes which are mutually parallel in all prime divisor
slices $\CS_D$;
see the figure shown in the introduction.
However, we do not exclude the case of $u^{ij}=0$.
Then coherence means that $\Delta^i_D=\Delta^j_D$
whenever both polytopes are non-empty.

The divisorial fans we obtained in
Proposition \ref{toricppfan} are coherent. 
Here, we would like to raise the
opposite question:
If $\CS$ is a set of pp-divisors, then we denote by
$\langle\CS\rangle$ the set of all intersections of elements of $\CS$.
If $\CS$ is coherent, under which conditions 
does $\langle\CS\rangle$ become a divisorial fan?

\begin{proposition}
\label{prop-cohpp}
Let $\mathcal{S}=\{\mathcal{D}^i\}$ be a coherent set of
pp-divisors. Then $\langle\CS\rangle$ inherits the coherence 
as well as the pp-property.
Moreover, for any non-negative map
$ \mu \colon \{\mbox{\rm prime divisors on $Y$}\} \to 
\RR_{\geq 0} $,
% as in Definition \ref{def-wspc},
the slices 
$\langle\mathcal{S}\rangle_\mu$ are 
% not necessarily complete 
polyhedral subdivisions of $N_\QQ$.
% If, moreover, the separating vectors $u^{ij}$ 
% from the coherence definition
% are always non-zero, then 
% the full-dimensional cells carry unique labels.
\end{proposition}

\begin{proof}
First, it is easy to show that coherence survives under finite intersections:
Checking this comes down to considering elements
$\D^i,\D^j,\D^k\in\CS$ and trying to separate $\D^i\cap\D^j$ from
$\D^k$. If, say, $k=j$, then this is done by $u^{ij}$. On the other hand,
if $p\notin\{i,j\}$, then one takes $u^{ik}+u^{jk}$ instead.
\vspace{0.5ex}\\
Now, we will see why the pp-property remains valid for intersections
$\D^i\cap\D^j$ of $\CS$-elements.
Denoting $u:=u^{ij}$ from the definition of coherence
and $Z^{i/j}:=\bigcup_{\Delta_D^{i/j}=\emptyset} D$,
the loci of $\D^{i/j}$ are the semiprojective $Y^{i/j}:= Y\setminus Z^{i/j}$.
Since both $Z^i$ and $Z^j$ are supports of effective, semiample divisors,
we can use the sum of these divisors to show that $Z^i\cup Z^j$ is of the
same quality, i.e.\ $Y^i\cap Y^j\subseteq Y$ stays semiprojective.
Moreover, the locus $Y^{ij}$ of $\D^i\cap\D^j$ is
\[
Y^{ij} = (Y^i\cap Y^j) \setminus \bigcup_
{\hspace{-1em}\Delta_D^i\cap\Delta_D^j=\emptyset\hspace{-1em}} D
\renewcommand{\arraystretch}{1.3}
\begin{array}[t]{l}
\;=\; (Y^i\cap Y^j) \setminus \supp \sum_D (\min\langle \Delta^j_D,u\rangle
- \max\langle \Delta^i_D,u\rangle) D
\vspace{0.5ex}\\
\;=\; (Y^i\cap Y^j) \setminus \supp (\D^i(-u) + \D^j(u)).
\end{array}
\]
>From $-u\in\tail(\D^i)^\vee$ and $u\in\tail(\D^j)^\vee$, we obtain
that the divisors $\D^i(-u)$ and $\D^j(u)$ are semiample on
$Y^i$ and $Y^j$, respectively. Hence their sum is semiample on
$Y^i\cap Y^j$. Since $\D^i(-u) + \D^j(u)$ is also effective,
this shows the semiprojectivity of $Y^{ij}$. 
\vspace{0.5ex}\\
We may exploit another fact from this: The previous equation shows that
$\D^i(-u)=-\D^j(u)$ on $Y^{ij}$. Thus, on the locus of
$\D^i\cap\D^j$, the divisor $-\D^j(u)$ is semiample.
On the other hand, denoting $\D^j=\sum_D\Delta_D^j \otimes D$,
Lemma~\ref{lastlemma} tells us
that $\D':=\sum_D\face(\Delta_D^j,u)\otimes D$
leads to the evaluations
$\D'(u')=\D(u'+\ell u)|_{Y^{ij}} - \D(\ell u)|_{Y^{ij}}$ 
for $\ell\gg 0$.
Hence, they are semiample, too.
\vspace{0.5ex}\\
Eventually, we may assume $\CS=\langle\mathcal{S}\rangle$
to deal with the polyhedral subdivision
$\mathcal{S}_\mu$.
Under summation according to $\mu$,
the coherent separation of the polyhedra inside the coefficients
$\CS_P$ transfers to $\CS_\mu$ as
$\max\langle \Delta^i_\mu, u^{ij}\rangle < 
\min\langle \Delta^j_\mu, u^{ij}\rangle$
or
$\max\langle \Delta^i_\mu, u^{ij}\rangle=\min\langle \Delta^j_\mu,
u^{ij}\rangle=:c_\mu^{ij}$
with
$\,\Delta^i_\mu\cap [\langle \kbb, u^{ij}\rangle=c_\mu^{ij}]=
\Delta^j_\mu\cap [\langle \kbb, u^{ij}\rangle=c_\mu^{ij}]$.
% Now,
% if $\Delta^i_\mu$ and $\Delta^j_\mu$ were equal (but with the different labels
% $i$ and $j$),
% then being full-dimensional would imply that
% $\max\langle \Delta^i_\mu, u^{ij}\rangle > 
% \min\langle \Delta^j_\mu, u^{ij}\rangle$. This provides a contradiction.
\end{proof}

The main ingredient of divisorial fans is the face relation, i.e.\
the issue of open embeddings. It is one of the advantages of coherent pp-sets 
that the class of open embeddings among its elements
is easier to describe than in the general case.

\begin{lemma}
\label{lem-cohopen}
Let $\D$ be a pp-divisor on $Y$ with non-empty coefficients,
let $Z\subseteq Y$ be the support of an effective, semiample divisor.
For $u\in\tail(\D)^\vee$, we define
$\D'$ on $Y\setminus Z$ via the $u$-faces of the $\D$-coefficients.
If $\D'$ is assumed to be pp (i.e.\ if $-\D(u)$ is semiample on $Y\setminus Z$),
then $X(\D')\to X(\D)$ is an open embedding if and only if
\[
\forall y\in Y\setminus Z \hspace{0.5em}
\exists D_y\in |\D(\NN\cdot u)|:
y\notin \supp D_y =: Z_y\supseteq Z. 
\]
\end{lemma}

\begin{proof}
By definition of $\D'$, we know that 
$\Delta_y'=\face(\Delta_y,u)$, hence 
$\face(\Delta_y',u)=\face(\Delta_y,u)$
for all points $y\in Y\setminus Z$.
In particular, the condition in the lemma is stronger than
that of Proposition \ref{openemb} or Definition \ref{def:ppface},
hence it is sufficient for having an open embedding.
\vspace{0.5ex}\\
Out of necessity, let us assume that $\D'$ is a face of $\D$ in the sense of
Definition \ref{def:ppface}.
Then, for any $y \in \Loc(\mathcal{D}')=Y\setminus Z$, there are
$u_y \in \sigma^\vee \cap M$
%%%%%%%%%%%%%%%%%%%%
and $D_y\in |\D(u_y)|$ with
$y\notin \supp D_y\supseteq Z$ and
% $\face(\Delta_y,u)=
$\Delta_y'=\face(\Delta_y,u_y)$.
Since $u$ defines the same face of $\Delta_y$, both
$u_y$ and $u$ are contained in the interior of the normal cone
${\mathcal N}(\Delta_y',\Delta_y)$. Thus, we can find another
$u'\in\interior {\mathcal N}(\Delta_y',\Delta_y)$
such that $u_y+u'=k\cdot u$ for some $k\gg 0$.
The semiampleness of $\D(u')$ provides some
$E\in |\D(u')|$ avoiding $y$, and we may use
$D_y+E\in |\D(u_y)+\D(u')|= |\D(ku)|$ as the new divisor $D_y$.
\end{proof}

The openness condition of Lemma \ref{lem-cohopen} looks like 
asking for semiampleness of some divisor. This is indeed the case if
$Z=\emptyset$. In the general case, however, we can only formulate
a sufficient condition in these terms:
% \vspace{0.5ex}\\
%
The condition of Lemma \ref{lem-cohopen} is fulfilled whenever
there is an effective, semiample divisor $E$ with
$\supp E= Z$ and $\,(k\,\D(u)-E)$ being semiample for $k\gg 0$.
A very special example for this situation is when $E\sim\D(u)$ 
-- this is what happens in the case of the localizations described in
Proposition \ref{prop:localization}.
A second class of easy instances is, of course, when $Y$ is affine.
% then the semiampleness
% of Cartier divisors is no condition at all.
Summarizing our considerations so far, we obtain

\begin{corollary}
\label{cor-suffan}
Let $\mathcal{S}=\{\mathcal{D}^j\}$ be a coherent set of
pp-divisors such that, for any $i,j$, there is
an effective, semiample divisor $E^{ij}$ on $\Loc(\mathcal{D}^j)\subseteq Y$
with $\,\supp E^{ij} = \bigcup\{D\kst \Delta_D^i\cap\Delta_D^j = \emptyset\}$
% $Z^{vw}:=\bigcup\{D_i\kst \Delta_i^v\cap\Delta_i^w = \emptyset\}\subseteq Y$
and $k\,\D^j(u^{ij})- E^{ij}$ being semiample for $k\gg 0$.
Then $\langle\CS\rangle$ is a divisorial fan.
\end{corollary}

\begin{proof}
Corollary \ref{sec:cor-intersect-divisors} establishes the relation
between intersection of pp-divisors, i.e.\ of their polyhedral coefficients,
and the intersection of the corresponding affine $T$-varieties.
Hence, the claim follows from Proposition \ref{prop-cohpp} and, 
via the previous remarks, from
Lemma \ref{lem-cohopen}.
\end{proof}

We will conclude this section by taking a closer look at the 
situation of two special cases. If $Y$ is either a
smooth, projective curve (covering the one-codimensional torus actions)
or $Y=\PP^n$, then
everything becomes very clear -- and we 
even have a straight
characterization of the pp-ness of the elements of $\CS$:

% \begin{proposition}
% \label{prop-curvan}
% Let $\mathcal{S}=\{\mathcal{D}^j\}$ be a coherent set of
% polyhedral divisors on a smooth, projective curve $Y$ over $\KK$.
% Then, the elements $\D^j$ are pp-divisors
% if and only if
% $\,\Loc(\mathcal{D}^j)=Y$ implies
% $\deg\D^j\subsetneq \tail\,\D^j$ and, additionally,
% $\D^j(u)=0$ for
% any $u\in (\tail\,\D^j)^\vee$ with 
% $\min\langle \deg\D^j, u\rangle =0$.
% %
% Moreover, $\langle\CS\rangle$ becomes a divisorial fan 
% if and only if, additionally,
% % in case of $u=u^{ij}$, 
% $\,\Loc(\mathcal{D}^j)=Y$ and
% $\D^j(u^{ij})=0$ imply
% $\,\deg \D^i\cap\deg \D^j\neq \emptyset$.
% \end{proposition}

\begin{proposition}
\label{prop-curvan}
Let $\mathcal{S}=\{\mathcal{D}^j\}$ be a coherent set of
polyhedral divisors on a projective $Y$ with either
$\dim Y=1$ or $\,Y=\PP^n$.
Then, the elements $\D^j$ are pp-divisors
if and only if
$\,\deg\D^j\subsetneq \tail\,\D^j$ and 
{\rm (being automatically satisfied if $Y$ is rational)}
$\D^j(u)\sim 0$ for
any $u\in (\tail\,\D^j)^\vee$ with
$\,\min\langle \deg\D^j, u\rangle =0$.\\
Assuming this, $\langle\CS\rangle$ becomes a divisorial fan
if and only if
$\,\min\langle \deg\D^j, u^{ij}\rangle =0$ implies that
$\,\deg \D^i\cap\deg \D^j\neq \emptyset$.
\end{proposition}

\begin{proof}
First, in case of $\deg\D^j=\emptyset$, all conditions mentioned in
the proposition are automatically fulfilled. On the other hand, this case
means that $\,\Loc(\mathcal{D}^j)\subsetneq Y$, i.e.\ that
$\,\Loc(\mathcal{D}^j)$ is affine. In particular, all conditions characterizing
pp-properties or openness of maps among the corresponding affine $T$-varieties
are satisfied, too.\\
Thus, we may assume that $\deg\D^j\neq \emptyset$. 
In the proposition,
the condition on the single polyhedral divisors $\D^j$
translates into $\deg\D^j(u)>0$, i.e.\ ampleness, or
$\D^j(u)\sim 0$, where the latter cannot occur for
$u\in \interior(\tail\,\D^j)^\vee$. This does exactly characterize pp-ness.
Now, assuming that $\CS$ is a set of pp-divisors,
we see that, in the case of $\dim Y=1$ or $Y=\PP^n$, 
the condition from Lemma \ref{lem-cohopen} comes down to the ampleness
of $\D(u)$ or, alternatively, to $\D(u)\sim 0$ and $Z=\emptyset$. 
In particular,
the sufficient condition from Corollary \ref{cor-suffan}
is necessary, too.
Thus, to characterize the divisorial fan property, it remains to
ask for $\bigcup\{D\kst \Delta_D^i\cap\Delta_D^j = \emptyset\}=\emptyset$
whenever $\D^{j}(u^{ij})\sim 0$.
\end{proof}

\section{Separateness and Completeness}
\label{sec:sepcomp}

Separateness of a toric variety
is reflected by the fact that any two
cones of its fan admit a separating linear 
form (cutting out precisely their 
intersection).
Moreover, a toric variety is complete 
if and only if the cones of its fan cover 
the whole vector space.
In this section, we extend these two 
observations to the setting of divisorial fans.

The idea is to interprete the 
valuative criteria for separateness and 
completeness in our combinatorial terms.
We will work with divisorial fans on smooth
semiprojective varieties $Y$;
this is no loss of generality, because,
given a divisorial fan $\mathcal{S}$ on a 
singular $Y$, we may resolve singularities,
then pull back $\mathcal{S}$, and the 
new  divisorial fan defines the same $T$-variety.

Let us briefly fix the notation
concerning valuations.
As usual, we mean by 
a valuation of $\KK(Y) / \KK$,
where $Y$ is any variety,
a valuation 
$\mu \colon \KK(Y)^\ast \to \QQ$ 
of the function field
with $\mu = 0$ along $\KK$.
Moreover, we say that $y \in Y$  
is the (unique) center of $\mu$ if
the valuation ring 
$(\mathcal{O}_\mu,\mathfrak{m}_\mu)$ 
dominates the local ring 
$(\mathcal{O}_y,\mathfrak{m}_y)$,
this means that 
$\mathcal{O}_y \subseteq \mathcal{O}_\mu$ 
and 
$\mathfrak{m}_y = \mathfrak{m}_\mu \cap \mathcal{O}_y$
hold.

\begin{definition}
Let $Y$ be a variety, 
and let $\mu$ be a valuation
of $\KK(Y) / \KK$ with center
$y \in Y$.
Then there is a well defined 
group homomorphism
$$ 
\mu \colon \CDiv(Y) \ \to \QQ,
\quad
D \ \mapsto \ \mu(f), 
\quad
\text{where } D = \div(f) 
\text{ near } y 
\text{ with } f \in \KK(Y),
$$
and, for smooth $Y$, 
this provides a weight function
$ 
\mu \colon 
\{\text{prime divisors on } Y\}
\ \to \
\QQ.
$
\end{definition}

\begin{remark}
The homomorphism 
$\mu \colon \CDiv(Y) \to \QQ$
associated to a valuation $\mu$ 
of $\KK(Y)/Y$ 
% with center $y \in Y$
satisfies  
$\mu(D) \ge 0$
for all $D \ge 0$. 
\end{remark}

Recall from Section~\ref{sec:coherent} that, 
for a divisorial fan $\CS$,
% on a smooth semiprojective variety $Y$,
we have defined the notion of a slice.
For valuations $\mu$,
the directed systems $\mu(\CS)$ of polyhedra in
$N_\QQ$ are weighted versions of our former $\CS_y$.
In particular, they share the property of being a complex.

\begin{definition}
Let $Y$ be a smooth semiprojective variety,
and let $\mathcal{S}$ be a divisorial fan  on $Y$
with polyhedral coefficients living in $N_\QQ$.
\begin{enumerate}
\item 
We say that  $\mathcal{S}$ is 
{\em separated\/}, if for any two
$\mathcal{D}^i, \mathcal{D}^j \in \mathcal{S}$
and any valuation on $Y$,
we have
$ \,{\mu}(\mathcal{D}^i \cap \mathcal{D}^j)
= {\mu}(\mathcal{D}^i)
\ \cap \
{\mu}(\mathcal{D}^j)$.
\item
We say that $\mathcal{S}$ is 
{\em complete\/}, if
$Y$ is complete, $\mathcal{S}$ is separated,
and, for every valuation $\mu$, 
the slice $\mu(\mathcal{S})$ 
covers~$N_\QQ$.
\end{enumerate}
\end{definition}

% \begin{remark}
% Every coherent divisorial fan on
% a smooth semiprojective variety $Y$
% is separated.
% \end{remark} 

\goodbreak

\begin{remark}
Let $\mathcal{S}$ be a divisorial fan on
a smooth semiprojective variety $Y$.
\begin{enumerate}
\item
The divisorial fan $\mathcal{S}$ is separated if 
and only if for every valuation
$\mu$,
the slice $\mu(\mathcal{S})$ is a polyhedral 
subdivision.
\item
\correct{}{If $Y$ is a smooth curve or if
$\mathcal{S}$ is coherent, then $\mathcal{S}$ is automatically separated.}
\item
The divisorial fan $\mathcal{S}$ is complete if 
and only if $Y$ is complete and, for 
every valuation $\mu$,
the slice $\mu(\mathcal{S})$ is a complete 
polyhedral subdivision.
\item
\correct{}{If $Y$ is a smooth, complete curve, then $\mathcal{S}$ is complete
if all prime divisor 
slices of $\mathcal{S}$ cover the 
whole vector space.}
\end{enumerate} 
\end{remark}

% \begin{remark}
% \label{sec:rem-curve-complete-separated}
% For a divisorial fan $\mathcal{S}$ 
% on a smooth curve $Y$,
% every valuative slice is the slice
% of a multiple of a prime divisor. 
% So, the divisorial fan $\mathcal{S}$ is 
% automatically separated, and
% moreover, it is complete if $Y$
% is complete and all prime divisor 
% slices of $\mathcal{S}$ cover the 
% whole vector space. 
% \end{remark}

\begin{theorem}\label{thm:sepcompchar}
Let $Y$ be a smooth semiprojective variety,
let $\mathcal{S}$ be a divisorial fan  on $Y$,
and let $X$ be the associated prevariety.
\begin{enumerate}
\item
$X$ is separated if and only if $\mathcal{S}$ is
separated.
\item 
$X$ is complete if and only if $\mathcal{S}$
is complete.
\end{enumerate}
\end{theorem}

The proof of this result is based on 
a characterization of
existence of centers for valuations
of the function field of an affine 
$T$-variety in terms of its defining 
pp-divisor.
A first step is to understand the 
valuations themselves.

\begin{remark}
\label{center2grphom}
Let $\mathcal{D}$
be a pp-divisor with tail cone 
$\sigma \subseteq N_\QQ$ 
on a smooth semiprojective variety $Y$, 
and let $X$ be the associated affine 
$T$-variety.
Then $\KK(X)$ is the quotient field of
the Laurent polynomial algebra
\begin{eqnarray*}
\KK(Y)[M]
& = & 
\bigoplus_{u \in M}
\KK(Y) \mal 1^u,
\end{eqnarray*}
where, as usual, $M = \Hom(N,\ZZ)$ 
is the dual lattice.
Given a valuation $\mu$ of the function field
$\KK(Y)$
and a vector $v \in N$, we obtain a 
map 
$$
\nu_{\mu,v} \colon \KK(Y)[M] \ \to \ \QQ,
\qquad
{\textstyle \sum_i}\,  f_i \,1^{u_i}
\ \mapsto \ 
\min_i \big(\mu  (f_{i}) + \bangle{u_i,v}\big).
$$

% where $f_{u_i} = f_i 1^{u_i}$ 
% stands for a homogeneous 
% element of degree $u_i \in M$.
This map extends to a valuation 
of the field $\KK(X)$, and we have a
canonical injection 
$$ 
{\rm valuations}(\KK(Y)/\KK) \times N
\ \to \ 
{\rm valuations}(\KK(X)/\KK),
\qquad
(\mu,v)
\ \mapsto \
\nu_{\mu,v}.
$$
Conversely, any valuation $\nu$ on $\KK(X)/\KK$
coincides on the homogeneous elements 
of $\KK(Y)[M]$ with a unique $\nu_{\mu,v}$:
the data $\mu$ and $v$ are defined via
$$
\mu \ = \ \nu_{\vert \KK(Y)},
\qquad
\bangle{u,v} \ = \  \nu(1^u).
$$
Thus, on the homogeneous elements of 
$\KK(Y)[M]$, any valuation $\nu$ of 
$\KK(X)$ is uniquely represented by 
a valuation $\nu_{\mu,v}$;
we will denote this by 
$\nu \sim \nu_{\mu,v}$.
\end{remark}

\begin{lemma}
\label{lem:center}
Let $\mathcal{D}$ be a pp-divisor  
on a smooth semiprojective variety $Y$, 
let~$X$ be the associated affine 
$T$-variety, and consider a valuation 
$\nu \sim \nu_{\mu,v}$ on $\KK(X)/\KK$.
Then the following statments are equivalent.
\begin{enumerate}
\item
The valuation $\nu$ on $\KK(X)/\KK$ has a center $x \in X$.
\item
The valuation $\mu$ on $\KK(Y)/\KK$ has a center $y \in Y$
with $v \in {\mu}(\mathcal{D})$.
\end{enumerate}
\end{lemma}

\begin{proof}
As usual, let
$A = \Gamma(Y,\mathcal{A}) = \Gamma(X,\mathcal{O})$ 
denote the global ring
associated to the pp-divisor
$\mathcal{D}$.

Suppose that $\nu$ has a center $x \in X$.
Then we obtain $A \subseteq \mathcal{O}_\nu$,
which implies 
$A_0 \subseteq \mathcal{O}_\nu$.
Thus, there is a center in 
$Y_0 := \Spec(A_0)$
for the restriction of $\nu$ 
to $\KK(Y_0)$.
Since $\mu$ and $\nu$ coincide
on $\KK(Y_0)$ and $Y$ is projective 
over $Y_0$, the valuative criterion 
of properness~\cite[Thm~II.4.7]{hartshorne:ag} 
provides a center $y \in Y$ for
$\mu$.

In order to verify the desired property 
for the weight function $\mu$, 
suppose, to the contrary, that
$v \not\in {\mu}(\mathcal{D})$
holds.
Then, there is a linear form 
$u \in \sigma^\vee \cap M$,
where $\sigma$ stands for the tail cone
of $\mathcal{D}$,
such that  
$$ 
\bangle{u,v} 
\ < \ 
\min
\bangle{u,{\mu}(\mathcal{D})} 
\ = \
\mu(\mathcal{D}(u)).
$$
Replacing $u$ with a suitable positive 
multiple, we achieve that $\mathcal{D}(u)$ 
is Cartier and base point free on $Y$. 
Then we may choose a 
global section $f \in \Gamma(Y,\mathcal{D}(u))$
such that $f^{-1}$ is a local equation for
$\mathcal{D}(u)$ near $y$. We obtain
$$ 
\mu(f)
\ = \ 
-\mu(\mathcal{D}(u))
\ < \ 
-\bangle{u,v}.
$$
This implies $\nu_{\mu,v}(f) < 0$,
and thus 
$f \in \Gamma(Y,\mathcal{D}(u)) 
\subseteq \Gamma(X,\mathcal{O})$ 
does not belong to the valuation ring 
$\mathcal{O}_{\nu}$.
Consequently, $\mathcal{O}_{\nu}$ cannot 
dominate any local ring $\mathcal{O}_x$,
meaning that $\nu \sim \nu_{\mu,v}$ has no
center in $X$; a contradiction.
 
Now suppose that $\mu$ has a center 
$y \in Y$ and that  
$\mu \colon \CDiv(Y) \to \QQ$ is as in~(ii). 
Then $\mu(\D(u)) 
= 
\min\bangle{u,{\mu}(\D)} 
\le 
\bangle{u,v}$ 
holds for all $u \in \sigma^\vee \cap M$.
Thus, we have
\begin{eqnarray*}
f \; \in \; \Gamma(Y, \Of(\D(u)))
& \implies & \div(f) \; \ge \; -\D(u)
\\
& \implies & 
\mu(f) \; \ge \; \mu(-\D(u)) \; \ge \; \bangle{u,v}
\\
&\implies& 
\nu(f) \; = \; \mu(f) + \bangle{u,v} \ \ge \;  0
\\
&\implies& 
f \ \in \; \Of_\nu.
\end{eqnarray*}
This implies $A \subseteq \Of_\nu$. 
Thus, there is a prime ideal $\mathfrak{p} \subseteq A$
such that $A_\mathfrak{p}$ is dominated by $\Of_\nu$.
In other words: $\nu$ has a center $x \in X$.
\end{proof}

\begin{proof}[Proof of Proposition~\ref{thm:sepcompchar}]
We first treat separateness by applying the well-known
valuative criterion, 
cf.~\cite[Thm.~II.4.3]{hartshorne:ag}. 
It says that a prevariety $X$ 
is separated if and only if
every valuation of $\KK(X)/\KK$ 
admits at most one center in~$X$.

Let $\mathcal{S}$ be separated.
Consider a valuation $\nu \sim \nu_{\mu,v}$ 
of $\KK(X)/\KK$ with centers $x, x' \in X$. 
Then $x$ and $x'$ belong to affine charts 
$X(\mathcal{D})$ and $X(\mathcal{D}')$
with $\mathcal{D}, \mathcal{D}' \in \mathcal{S}$.
By Lemma~\ref{lem:center}, there
is a (unique) center $y \in Y$ for $\mu$,
and we have 
\begin{eqnarray*}
v
& \in & 
{\mu}(\mathcal{D}) 
\ \cap \ 
{\mu}(\mathcal{D}'). 
\end{eqnarray*}
Separateness of $\mathcal{S}$ gives 
$v \in {\mu}(\mathcal{D} \cap \mathcal{D}')$.
Again by Lemma~\ref{lem:center},
it follows that $\nu$ has a center in 
$X(\mathcal{D} \cap \mathcal{D}') 
= X(\mathcal{D}) \cap X(\mathcal{D}')$.
Since $X(\mathcal{D})$ and $X(\mathcal{D}')$
are separated, this implies $x = x'$.

Now, let $X$ be separated,
and suppose that $\mathcal{S}$ 
is not.
Then there are 
$\mathcal{D}, \mathcal{D}' \in \mathcal{S}$
and a valuation $\mu$ on $Y$ such that we have
\begin{eqnarray*}
{\mu}(\mathcal{D} \cap \mathcal{D}')
& \subsetneq & 
{\mu}(\mathcal{D})
\ \cap \
{\mu}(\mathcal{D}').
\end{eqnarray*}
Pick any 
$v \in {\mu}(\mathcal{D}) 
\cap {\mu}(\mathcal{D}')$
 that does not belong to 
${\mu}(\mathcal{D} \cap \mathcal{D}')$.
Then by Lemma~\ref{lem:center},
the valuation $\nu = \nu_{(\mu,v)}$ 
has no center in $\ppv{\D \cap \D'}$ but 
it has centers $x \in \ppv{\D^\prime}$ 
and $x' \in \ppv{\D^{\prime\prime}}$.
Since $\ppv{\D \cap \D'}$ equals 
$\ppv{\D }\cap \ppv{\D'}$,
we obtain $x \ne x'$.
This contradicts separateness of $X$.

We turn to completeness.
Again we make use of a valuative criterion,
cf.~\cite[Thm.~II.4.7]{hartshorne:ag}.
It says that a (separated) variety $X$ is 
complete if and only if every valuation of
$\KK(X)/\KK$ has a center in $X$.

Let $X$ be complete.
Suppose that $Y$ is not. 
Then there is a valuation $\mu$ of $\KK(Y)/\KK$ 
without center on $Y$.
Take any $v \in N_\QQ$. 
Lemma~\ref{lem:center} says that
$\nu_{\mu,v}$  
has no center in $X$.
A contradiction.
Next suppose that there were a 
noncomplete valuative slice $\mu(\CS)$.
Take $v \in  N_\QQ \setminus |\mu(\CS)|$. 
Then, again by Lemma~\ref{lem:center},
the valuation $\nu_{\mu,v}$ on $\KK(X)/\KK$ 
has no center in $X$.
A contradiction.

Finally, let $\mathcal{S}$ be complete.
Let $\nu \sim \nu_{\mu,v}$ be any valuation 
on $\KK(X)/\KK$.
By completeness of $Y$, the valuation 
$\mu$ has a center
$y \in Y$.
Moreover, the slice $\mu(\CS)$ is complete, 
this means that there is a $\D \in \CS$
such that $v \in {\mu}(\D)$ holds.
Lemma~\ref{lem:center} then provides 
a center of $\nu$ in $\ppv{\D} \subseteq X$. 
\end{proof}

\begin{remark}
Example \ref{ex-noncoh} with, e.g.,
$Y=\AA^1$ and $D_1$, $D_2$ being two different prime divisors,
yields a separated, but non-coherent divisorial fan.
\end{remark}

The following two examples underline
that for checking separateness and 
completeness it is not sufficient
to consider the prime divisor slices.

\begin{example}
Let $Y := \mathbb{A}^2$ with the 
coordinate functions $z$ and $w$, 
and denote the coordinate axes by
$D := \div(z)$ and $E := \div(w)$.
We consider pp-divisors
$$ 
\mathcal{D}^1
\ = \ 
\Delta^1_D \otimes D
\; + \; 
\Delta^1_E \otimes E,
\qquad
\mathcal{D}^2
\ = \ 
\Delta_D^2 \otimes D
\; + \; 
\Delta_E^2 \otimes E
$$
with tail cones 
$\tail(\mathcal{D}^1) = \QQ_{\le 0} \mal e_1$ 
and 
$\tail(\mathcal{D}^2) = \QQ_{\ge 0} \mal e_1$
and polyhedral coefficients in $\QQ^2$
according to the following figure:

\bigskip

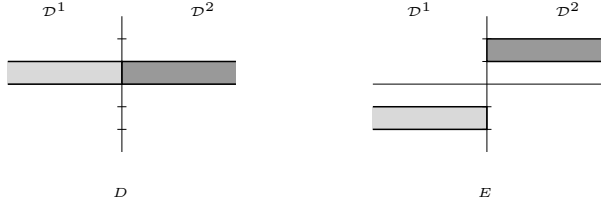
\begin{center}
\input{nonsep2.pstex_t}
\end{center}
Then $\mathcal{S} := \{\mathcal{D}^1, \mathcal{D}^2\}$
generates a divisorial fan, and the above figure shows 
that its prime divisor slices 
are nice polyhedral subdivisions.
However, the divisorial fan $\mathcal{S}$ is not 
separated.
While the weight function $\mu_{(0,0)}$ yields the subdivision
$\CS_{(0,0)}$ consisting of the polytopes
$\D^i_{(0,0)}=\Delta^i_D+\Delta_E^i$ and $\emptyset$,
the valuation 
$$
\textstyle
\mu(\sum_{a,b} \lambda_{a,b} z^a w^b):=
\min\{2a+b\kst \lambda_{a,b} \neq 0\}
$$
provides $\mu(\CS)=\{\mu(\D^1),\mu(\D^2),\mu(\D^1\cap\D^2)\}$  
with
$$ 
\mu(\mathcal{D}^1 \cap \mathcal{D}^2)
\ = \ 
\emptyset,
\hspace{0.5em}\mbox{but}\hspace{0.8em}
\mu(\mathcal{D}^1) \;  \cap  \; \mu(\mathcal{D}^2)
\ = \ 
\{(0,1)\}.
$$
One also sees directly that 
$X(\mathcal{S})$, as the gluing
of $\Spec(\KK[x,y,s,ty^{-1}, xy^2t^{-1}])$ 
and $\Spec(\KK[x,y,s^{-1},ty^2, xt^{-1}y^{-1}])$
along $\Spec(\KK[x,y,y^{-1},s,s^{-1},t,xt^{-1}])$,
is not separated.
\end{example}

\begin{example}
On $Y := \PP^2$, fix two coordinate axes, 
say $D$ and $E$, and consider the coherent set 
$\mathcal{S} = \{\mathcal{D}^1, \ldots, \mathcal{D}^4\}$ of polyhedral
divisors
given by its prime divisor slices as indicated below.
By Proposition \ref{prop-curvan}, we know that
$\bangle{\CS}$ is a divisorial fan. 
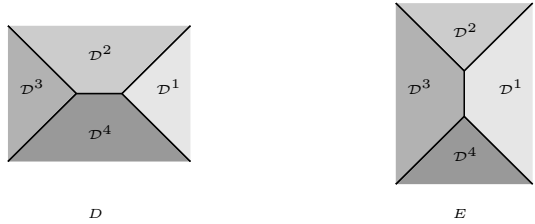
\begin{center}
\input{noncompl.pstex_t}
\end{center}
The prime divisor slices are complete, but there is
a valuation $\mu$ with $\mu(D)=\mu(E)=1$
providing the following
noncomplete 
slice.
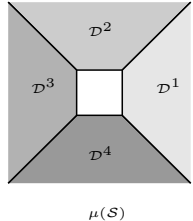
\begin{center}
\input{noncompl2.pstex_t}
\end{center}
%\kq{Die Beschriftung mu\ss\ zu $\mu(\CS)$ geaendert werden.}\\
%But this fan can be compactified, by adding the polyhedral divisor
%$$\D^\infty \quad = \quad- \otimes D \quad + \quad | \otimes E \quad + \quad \emptyset \otimes F.$$
\end{example}

\section{Further examples}
\begin{example}
\label{ex-flenner}
In \cite[Example 5.13]{Flenner:complete}, Flenner et al.\ describe
the so-called Danilov-Gizatullin $\KK^\ast$-surfaces $V_0$.
They depend on the choice of $r,s\in\ZZ_{\geq 1}$ and points
$p_0,p_1\in\KK^1$, and, using our language of pp-divisors, are
given by
$$
\D^0=[-\frac{1}{r},0]\otimes \{p_0\} + [0,\frac{1}{s}]\otimes \{p_1\}
\in \PPDiv_\QQ(\KK^1,\{0\})
$$
with $N=\ZZ$. In \cite{Flenner:complete},
the surface $V_0$ is described by the two
$\QQ$-divisors $D_+=\D^0(1)=-\frac{1}{r} p_0$ and
$D_-=\D^0(-1)=-\frac{1}{s} p_1$ on $\KK^1$.
Proposition 3.8 of \cite{Flenner:complete} provides 
a canonical $\KK^\ast$-equivariant completion of $V_0$ by ``adding''
two charts $V_+$ and $V_-$. Using our divisorial fans,
this completion is given as $\widetilde{X}(\CS)$
% of the following divisorial fan 
with $\CS=\{\D^-,\D^0,\D^+\}$
\begin{center}
\input{gizatul.pstex_t}
\end{center}
(and $\D^0_\infty=\emptyset$) on $\PP^1_\KK$.
Note that the contraction $X(\CS)$ of 
$\widetilde{X}(\CS)$ provides an even
``smaller'' compactification of $V_0$.
The completeness of $X(\CS)$ and $\widetilde{X}(\CS)$ is reflected by the fact
that, for every prime divisor, 
$M=\QQ$ is completely covered by the corresponding polyhedral coefficients 
(cf.\ Theorem~\ref{thm:sepcompchar}).
% and Remark~\ref{sec:rem-curve-complete-separated}).
\end{example}

\begin{example}
\label{ex-vollmert}
In \cite[IV.1]{kempf:toroidal} Mumford constructed a toroidal structure
on normal varieties $X$ with one-codimensional torus action. 
% In particular, associated to those objects,
% there is an abstract fan carrying some information about $X$.
This leads to an associated fan carrying some information about the
original $X$.
However, in contrast to the case of toric varieties,
this fan is not sufficient to recover $X$.
 
Since the theory of divisorial fans developed in the present paper 
provides a tool for keeping complete information about the normal varieties
with their torus action, there should be a direct way to translate
divisorial fans on curves $Y$ into the Mumford fans. 
For affine $X$, this has been done by Vollmert in \cite{Robert}:
The Mumford fan is obtained by gluing the homogenizations of the polyhedral
coefficients of the pp-divisor along the tail cone 
that is their common face.
As Mumford's construction works for non-affine varieties, too,
one easily sees that Vollmert's theorem remains valid in the general case.
\end{example}

\label{sec:txt-bundles}
In the next examples, 
we consider equivariant 
vector bundles on a given 
toric variety
$X = \toric{\Sigma}$ 
arising from a fan $\Sigma$ 
in $N_\QQ$.

Recall that Klyachko~\cite{klyachko:eqbundles}, 
and later Perling~\cite{perling:sheaves}, 
gave a combinatorial description 
of the $T$-equivariant, reflexive 
sheaves $\CE$ on $X$.
In Perling's notation,  $\CE$ 
is given by a $\KK$-vector 
space $E$ together
with $\ZZ$-labeled increasing filtations 
$E^\varrho(i)$ for every $\varrho\in\Sigma^{(1)}$. 
Note that not only is the filtrations itself 
an important data, but also the $i$, 
which tells you when a jump takes place.

The sheaf $\CE$ is locally free 
(of rank $ r =\dim_k E$) if and only if 
we can find for every $\sigma \in \Sigma$ 
a basis $e^\sigma_1 \ldots e^\sigma_r$ of 
$E$ and weights 
$u^\sigma_1, \ldots, u^\sigma_r \in M$ such that 
for all $\varrho \in \sigma(1)$ one has
\begin{eqnarray*}
e^\sigma_j \in E^\varrho(i) 
& \iff &
\mul{u^\sigma_j}{\varrho} \geq i.
\end{eqnarray*}
For the affine charts $\toric{\sigma} \subseteq X$,
we have 
$\CE(\toric{\sigma}) \subseteq \CE(T) \cong  E \otimes_\KK \KK[M]$; 
thus these data corresponds directly to 
elements of the graded module $\CE(\toric{\sigma})$ 
which generate $\CE|_{\toric{\sigma}}$ freely.

There are two striking examples for this
encoding via filtrations.
First, reflexive sheaves of rank one utilize 
a one-dimensional vector space~$E$. 
Since their filtrations are completely determined by telling
for which $i_\varrho$ the unique step 
$$
E^\varrho(i_\varrho-1)=0, 
\qquad
E^\varrho(i_\varrho)=E$$ 
takes place, 
Klyachko's description just means to fix a map $i:\Sigma^{(1)}\to\ZZ$.
This coincides with the classical description, see, e.g.,
\cite{kempf:toroidal}. 
Second, the cotangent bundle $\Omega_{X}$ on a smooth 
$X = \toric{\Sigma}$ may be obtained from the vector space 
$E:=M \otimes_\ZZ \KK$ with the filtrations
$$
E^\varrho(0)=E,
\qquad
E^\varrho(1) = \KK \cdot \varrho^\perp, 
\qquad
E^\varrho(2) = E.
$$

%%!!! nochmal pruefen: $i$ oder $-i$, $\varrho$ oder $\varrho^\perp$ !!!

%%%%%%%%
%  tact
%%%%%%%%%%
If $\CE=\oplus_j\mathcal{L}_j$ is a splitting vector bundle on 
$X=\toric{\Sigma}$, then
$\PP(\CE)$ is toric again
(under a larger torus). Assume, e.g., that the direct summands 
$\mathcal{L}_j$ of $\CE$
are ample and given by lattice polytopes $\Delta_j\subseteq M_\QQ$.
%If $\otimes_j\mathcal{L}_j=\Of(\sum_j \Delta_j)$ is ample on $X$, 
Then, $\PP(\CE)$ is associated to the normal fan
of $$\tilde{\Delta}:=\conv(\bigcup_j \Delta_j\times\{e_j\})\subseteq
M_\QQ\times\QQ^m.$$ 
%This is related to the so-called Cayley trick and may be found in (???).

If $\CE$ does not split, then $\PP(\CE)$ admits a torus
action of our original, lower-dimensional $T$. This suggests understanding
$\PP(\CE)$ as a divisorial fan $\CS$ on some variety $Y$. 
The slices of $\CS$ should be polyhedral subdivisions of $N_\QQ$.

\begin{example}
Consider an equivariant rank 2 bundle $\CE$ on a toric variety, 
given by some vector space filtrations $E^\varrho(i)$ of $ E = \KK^2$. 
Set $Y = \PP^1 \cong \PP(E^*)$. 
For every maximal cone $\sigma$, 
we obtain coefficients for two 
polyhedral divisors by cutting $\sigma$ 
with affine hyperplanes orthogonal 
to $u^\sigma_1 - u^\sigma_2$.
$$
  \Delta_\sigma^{1} =  \{ v \in N_\QQ \mid \mul{u^\sigma_1 - u^\sigma_2}{v} \geq 1 \} \cap \sigma , \ 
  \Delta_\sigma^{2} =  \{ v \in N_\QQ \mid \mul{u^\sigma_2 - u^\sigma_1}{v} \geq 1 \} \cap \sigma $$
$$
  \nabla_\sigma^{1} =  \{ v \in N_\QQ \mid \mul{u^\sigma_1 - u^\sigma_2}{v}  \leq 1 \} \cap \sigma, \ 
  \nabla_\sigma^{2} =  \{ v \in N_\QQ \mid \mul{u^\sigma_2 - u^\sigma_1}{v}  \leq 1 \} \cap \sigma
$$
The participating prime divisors for the $\sigma$-chart are $\{(e_1^\sigma)^\perp\}, \{(e_2^\sigma)^\perp\} \in \PP^1$, where $(e_i^\sigma)^\perp$ is the one dimensional subspace of $E^*$ orthogonal to $e_i$ and thus a element of $Y$. 
To be precise, we define the polyhedral divisors $\D^+_\sigma,\D^-_\sigma$ on $Y=\PP^1$.
\begin{eqnarray*}
  \D^+_\sigma & = & \Delta^{1}_\sigma \otimes \{(e_{1}^\sigma)^\perp \} + \Delta^{2}_\sigma \otimes \{(e^\sigma_{2})^\perp\}\\
  \D^-_\sigma & = & \nabla^{1}_\sigma \otimes \{(e_{1}^\sigma)^\perp \} + \nabla^{2}_\sigma \otimes \{(e_{2}^\sigma)^\perp\}.
\end{eqnarray*}

\begin{proposition}
  The set $\{\D^\pm_\sigma \mid \sigma \in \Sigma^{\text{max}}\}$ of all these polyhedral divisors generates a fan which encodes $\PP(\CE)$.
\end{proposition}
\begin{proof}
  We consider a maximal cone $\sigma \in \Sigma^{\text{max}}$. For now we fix the coordinates of $\PP^1$ given by the dual basis $(e^\sigma_1)^*,(e^\sigma_2)^*$. Then we have
  $$\D^+_\sigma  =  \Delta^{1}_\sigma \otimes \{ 0 \}  + \Delta^{2}_\sigma \otimes \{ \infty \}.$$

  We consider $\KK[M][\frac{x_1}{x_2},\frac{x_2}{x_1}]$ with the obvious $M$-grading. Then the elements with weight $u$ are exactly those of the form $\chi^{u}f$ with $f \in \KK[\frac{x_1}{x_2},\frac{x_2}{x_1}] \subset \KK(\PP^1)$. With this representation, we have the graded isomorphisms
  \begin{eqnarray*}
    \KK[\sigmav \cap M][\chi^{u^\sigma_1 - u^\sigma_2} \frac{x_1}{x_2}] &\rightarrow& \bigoplus_u \Gamma(\Of(\D^+_\sigma(u))) \\
    \chi^{u}f &\mapsto& f \in \Gamma(\Of(\D^+_\sigma(u)))
  \end{eqnarray*}
    
    \begin{eqnarray*}
      \KK[\sigmav \cap M][\chi^{u^\sigma_2 - u^\sigma_1} \frac{x_2}{x_1}] &\rightarrow& \bigoplus_u \Gamma(\Of(\D^-_\sigma(u))) \\
      \chi^{u}f &\mapsto& f \in \Gamma(\Of(\D^-_\sigma(u)))
    \end{eqnarray*}

    \begin{eqnarray*}
      \KK[\sigmav \cap M][\chi^{u^\sigma_1 - u^\sigma_2} \frac{x_1}{x_2}, \chi^{u^\sigma_2 - u^\sigma_1} \frac{x_2}{x_1}] &\rightarrow& \bigoplus_u \Gamma(\Of((\D^+_\sigma \cap \D^-_\sigma)(u))) \\
      \chi^{u}f &\mapsto& f \in  \Gamma(\Of((\D^+_\sigma \cap \D^-_\sigma)(u)))
    \end{eqnarray*}
    
    Thus, we obtain $\PP(\CE|_{X_{\sigma}})$ by gluing $\ppv{\D^+_\sigma}$
    and $\ppv{\D^-_\sigma}$. For getting the global result we also need to take into account the base change $e^\sigma_1, e^\sigma_2 \mapsto e^\delta_1, e^\delta_2$ between two cones.
\end{proof}  
\end{example}

\begin{example}[Cotangent bundle]
%
%For equivariant bundles of higher rank $Y$ becomes a blowing up of $\PP^n$ and the pp-divisors get slightly more complicated.
%But for the cotangent bundles no blowing up is needed, so we can give a nice decription by divisorial fans.
%
Let $X=\toric{\Sigma}$ be a smooth toric variety. 
For every (maximal) cone $\sigma \in \Sigma$ 
we consider the polyhedra:
$$
\Delta_\sigma^{ij} 
= 
\bigcap_{k = 1}^n\{ v \in N_\QQ \mid \bangle{u^\sigma_k - u^\sigma_i, v} 
\geq 
\delta_{ij} - \delta_{jk}\} \cap \sigma,
$$
where $u^\sigma_1, \ldots u^\sigma_n$ 
is a $\ZZ$-basis of $M$, 
such that the dual basis  $(u^\sigma_1)^* \ldots (u^\sigma_n)^*$ 
contains the primitive generators 
of the rays in $\sigma(1)$. 
We set $Y = \PP(N \otimes \KK)$ and 
$$
\D_\sigma^{\,i} 
= 
\sum^n_{j=1}  \Delta_\sigma^{ij} \otimes (u^\sigma_j)^\perp.
$$
Then the divisorial fan $\CS_\Omega(\Sigma)$ 
generated by the set $\{\D_\sigma^{\,i}\}_{\sigma,i}$ 
describes the cotangent bundle on $X$. 
Note that the tail fan of $\CS_\Omega(\Sigma)$ 
is obtained by a barycentric subdivision 
of the $n$-dimensional cones of $\Sigma$.

Here are two concrete examples.
Firstly, for $\Omega_{\PP^2}$ on $\PP^2$ we 
obtain the picture already shown 
in the introduction.
% \bigskip
% \begin{center}
% \input{tangp2.pstex_t}
% \end{center}
% \\[1.0ex]
Secondly, for the cotangent bundle $\Omega_{\text{dP}_6}$ 
on the del Pezzo surface $\text{dP}_6$,
we obtain the divisorial fan given 
by following prime divisor slices.
\begin{center}
\bigskip
\input{delpezzo.pstex_t}
\end{center}
\end{example}

\goodbreak

\end{document}